\documentclass[11pt,a4paper]{elsarticle}
\usepackage{lipsum}
\makeatletter
\def\ps@pprintTitle{%
 \let\@oddhead\@empty
 \let\@evenhead\@empty
 \def\@oddfoot{}%
 \let\@evenfoot\@oddfoot}
\makeatother

\usepackage{graphicx}
\usepackage{graphics}
\usepackage{amsmath}
\usepackage{times}
\usepackage{psfrag}
\usepackage{amssymb}
\usepackage{subfigure}
\usepackage{bm}
\usepackage{color}
\usepackage[dvipsnames]{xcolor}

\begin{document}
	
	\begin{frontmatter}
		
		\title{Modeling of thin plate flexural vibrations by Partition of Unity Finite Element Method}
		
		\author[PolyU]{Tong~Zhou}
		
		\author[UTC]{Jean-Daniel~Chazot}
		
		\author[UTC]{Emmanuel~Perrey-Debain}
		
		\author[PolyU]{Li~Cheng\corref{cor1}}
		\ead{li.cheng@polyu.edu.hk}

		\address[PolyU]{Department of Mechanical Engineering, The Hong Kong Polytechnic University, Hong Kong, China}
		\address[UTC]{Laboratoire Roberval, FRE UTC-CNRS 2012, Universit\'e de Technologie de Compi\`egne,\\60205 Compi\`egne, BP 20529, France}

		\cortext[cor1]{Corresponding authors}
		
		\begin{abstract}
			
	This paper presents a conforming thin plate bending element based on the Partition of Unity Finite Element Method (PUFEM), for the simulation of steady-state forced vibration.	The issue of ensuring the continuity of displacement and slope between elements is addressed by the use of cubic Hermite-type Partition of Unity (PU) functions. With appropriate PU functions, the PUFEM allows the incorporation of the special enrichment functions into the finite elements to better cope with plate oscillations in a broad frequency band.	The enrichment strategies consist of the sum of a power series up to a given order and a combination of progressive flexural wave solutions with polynomials.	The applicability and the effectiveness of the PUFEM plate elements is first verified via the structural frequency response.	Investigation is then carried out to analyze the role of polynomial enrichment orders and enriched plane wave distributions for achieving good computational performance in terms of accuracy and data reduction.	Numerical results show that the PUFEM with high-order polynomials and hybrid wave-polynomial combinations can provide highly accurate prediction results by using reduced degrees of freedom and improved rate of convergence, as compared with the classical FEM. 
	
	
	

		\end{abstract}

		
	\end{frontmatter}

	\section{Introduction}


	Finite Element Method (FEM) is an important numerical simulation tool that is widely used for analyses of acoustics and structural vibrations.
	Mainly used for low frequency applications, it shows its limitations when dealing with short-wave simulations which involve a large number of wavelengths in the problem domain, thus leading to exorbitant computational cost.
	This results from the rule of thumb that a certain fixed number (around 10) of nodal points is required to model each wavelength.
	Moreover, a refined mesh discretization is also needed to control the pollution error and maintain good numerical accuracy, especially when the wave number becomes very high.

    The last decade has witnessed the rapid development of the Partition of Unity Finite Element Method (PUFEM) \cite{Melenk, BABUSKA1997} in solving the problem of short-wave modeling.
	This enriched method relies mainly on revamping the conventional finite element shape functions by using special auxiliary functions that can enhance the element approximation capability for the solution of the considered problem.
	The PUFEM has been successfully applied to simulate the acoustic field and structural response with short wavelength and strong oscillations \cite{Chazot2013_pufem1, Chazot2014_pufem2, YANG2018, ElKacimi2009, zhou2019, zhou2020}.
	It is shown that the PUFEM element can offer highly accurate numerical results with significantly improved convergence rate and a huge reduction in the number of degrees of freedom, as compared to the classical FEM.

	The objective of the present work is to investigate the applicability and the effectiveness of the PUFEM for solving the flexural vibration response of thin plate structures.
	The modeling and analysis of such structures covering a broad frequency range is an issue of growing importance for vibration and noise control applications \cite{cuenca2012, zhang2017, rouch2003, vanmaele2007, geng2017}.
	Traditionally, the vibration analysis of a plate structure is usually carried out by employing the thin plate bending elements based on the Kirchhoff-Love theory \cite{petyt2010, zienkiewicz2013}.
	The bending deformation of the thin plate is described by the lateral displacement of the plate middle surface.
	It is necessary to ensure the continuity of displacement and slope between adjacent elements ($C^1$ continuity) since the highest derivatives in the weak variational form for thin plate bending are second order \cite{strang1973, zienkiewicz2013}.
	The elements satisfying this requirement are referred to as conforming elements.
	One most commonly-used conforming classical thin plate element is constructed based on Hermite polynomial shape functions \cite{bogner1965, bardell1991, beslin1997}. 
	There are few works on the development of high-order shape functions with $C^1$ continuity compared with those of $C^0$ family \cite{zienkiewicz2013}, which only need displacement continuity. This stems from the additional continuity restriction for the slope \cite{strang1973, hughes1977, petyt2010}.
	Numerous efforts have been made in the literature to find alternative approaches to apply FEM framework for thin plate bending \cite{melosh1963, deak1967, batoz1982, petyt2010, zienkiewicz2013}, but few of them can be easily extended to achieve high-order approximation.

	The PUFEM technique investigated in this work offers a flexible and versatile platform for the design and development of new elements for many wave problems \cite{Bettess}.
	One essential component of the PUFEM element is the Partition of Unity (PU) function \cite{BABUSKA1997, babuvska2004}, which determines the smoothness of the element shape functions \cite{babuvska2004}.
	Most PUFEM works on acoustics \cite{Chazot2013_pufem1, YANG2018, christodoulou2017} and structural vibrations \cite{HAZARD2007, DeBel2005, zhou2019} adopt the piecewise Lagrangian FEM shape functions as PU functions. In most cases the linear shape functions are chosen. These functions only guarantee $C^0$ continuity, since a jump of the slope occurs at element interfaces. 
	The Hermite shape functions used in classical FEM satisfy the displacement and slope continuity requirement, but these functions do not form partition of unity \cite{zhou2019}.
	However, the components of Hermite functions only associated with the nodal displacement form a partition of unity and we can use these $C^1$ PU functions for developing PUFEM thin plate elements \cite{babuvska2004}.
	While the issues of constructing smooth basis functions with higher order continuity for enriched elements have been addressed by several previous works \cite{babuvska2004, oh2008, oh2012, davis2014, de2009ac}, research on the PUFEM with cubic Hermite-type PU functions for modeling flexural wave propagation in a thin plate is still lacking.

	
	
	

	In this article, we focus on the PUFEM formulation of a conforming thin plate bending element for steady-state forced vibration analyses.
	The continuity of the displacement and the slope between elements is imposed by the Hermite-type $C^1$ PU functions.
	With the appropriately enforced inter-element continuity, the PUFEM plate element allows the incorporation of auxiliary enrichment functions with good approximation properties for simulating broadband plate vibrational response.
	The enrichment strategies using polynomial sequences up to a given order and a combination of free wave solutions to the plate bending equation with additional polynomials are employed in order to better cope with the oscillatory behavior of the plate vibration.
	The performance of different enriched elements is evaluated in terms of accuracy and data reduction.
	The proposed PUFEM formulation is first applied for solving one-dimensional thin plate bending problem, where the flexural motion depends only on one spatial coordinate and is independent of other directions.
	A two-dimensional rectangular PUFEM plate element is then constructed for the general plate vibration analysis.
	In such a case, the flexural waves can travel in any direction over the planar plate and the wave propagation direction becomes an important factor governing the structural dynamic behaviors. 
	Particular attention is devoted to the development and investigation of a hybrid wave-polynomial enrichment that combines a set of plane-progressive flexural waves in various directions over a plane and additional polynomial terms.
	Increasing the plane wave distributions attached at each node leads to a hierarchic refinement for bidimensional wave enrichment field. 
	All the plate boundary conditions are set to be simply supported so the numerical results can be checked against the analytical solutions by the modal superposition method.
	Classical FEM using Hermite shape functions is also studied for comparison purposes.
	We adopt the Lagrange multiplier technique to impose the boundary conditions for all the elements in this work.	


	
	The remaining part of the paper is organised in the following way.
	In Section 2, the thin plate bending problem is solved by the variational approach in 1D configuration.
	This simple unidimensional case permits easy demonstration on the necessity of imposing the $C^1$ continuity by the Hermite-type PU functions.
	The performance of pure polynomial enrichment of different orders is studied and the advantage of using hybrid enrichment that contains a pair of propagating waves in opposite directions is illustrated without melting waves in other directions.
	Section 3 applies and extends the PUFEM formulation to obtain two-dimensional solutions to thin plate bending vibrations. 
	The PU functions and enrichment strategies are adapted to the case with two coordinate variables and a discretization scheme for Lagrange multipliers along plate borders is also given.
	Section 4 presents the numerical results calculated by two-dimensional PUFEM elements with hybrid enrichment.
	This section demonstrates the effectiveness of hierarchic refinement for wave field and the contributions of the added polynomials for plate displacement approximation.
	The influence of mesh irregularity is also under investigation.
	Conclusions are drawn in the final section.
	

	\section{One-dimensional solution}
	Consider the flexural vibration of a thin flat plate with thickness $H$, which is small compared with the dimensions defining the span of the surface lying in  the $(x,y)$ plane.
	The plate is assumed to be excited by an out-of-plane harmonic load at a circular frequency $\omega$ and the time factor $e^{-j \omega t}$ is implicit hereafter ($j=\sqrt{-1}$). 
	In the thin plate theory, the normals to the mid-surface of the undeformed plate are assumed to remain normal to it during the deformation process and the plate deformation can be described only by the lateral displacement $W$ of the structural mid-surface located at $z=0$. In this section we are interested in one-dimensional solutions which do not depend on the $y$ direction. In this case, the equation of motion for the plate has the same form as the Euler-Bernoulli beam equation:
	\begin{equation}\label{EoM_1d}
	D  \frac{\partial ^4 W}{\partial x^4}  - \omega^2 \rho H  W = F \delta (x-x_F),
	\end{equation}
	where $D=EH^3/12(1-\nu^2)$ is the bending rigidity of the plate per unit width, with Young's modulus $E$ and Poisson's ratio $\nu$, and $\rho$ is the material density. A concentrated force of magnitude $F$ at $x_F$ is expressed by using the spatial Dirac delta function $\delta(\cdot)$.
	With simply supported boundary conditions, the lateral displacement and bending moment vanish at both extremities, i.e. at $x=0,L_x$. In this scenario, the constrained variational formulation becomes
	\begin{equation}\label{Variation_1d}
	\int  _0^{L_x} \left( \frac{\partial ^2 \delta \hspace{-0.2mm} W}{\partial x^2} D  \frac{\partial ^2 W}{\partial x^2} - \omega^2 \delta \hspace{-0.2mm} W \rho H W    \right) \mathrm{d}x -\delta \hspace{-0.2mm} W|_{x_F} F - \delta \hspace{-0.2mm} W|_0 \mathit{\Lambda}|_0 - \delta \hspace{-0.2mm} W_{L_x} \mathit{\Lambda}|_{L_x}  = 0,
	\end{equation}
	where $\delta W$ is the arbitrary virtual displacement and the subscript of $W$ represents the position along the $x$ coordinate. In the above formulation, the boundary shear forces appear naturally as Lagrange multipliers $\mathit{\Lambda}|_0$ and $\mathit{\Lambda}|_{L_x}$. The essential boundary conditions can be weakly enforced as \cite{Laghrouche2, Chazot2013_pufem1, zhou2019}: 
	\begin{eqnarray}\label{BC_1d}
	\delta \hspace{-0.2mm} \mathit{\Lambda}|_0 W|_0 = \delta \hspace{-0.2mm} \mathit{\Lambda}|_{L_x} W|_{L_x}  =0, \hspace{0.6cm}  \forall (\delta \hspace{-0.2mm} \mathit{\Lambda}|_0, \delta \hspace{-0.2mm} \mathit{\Lambda}|_{L_x}).
	\end{eqnarray}
	

	\subsection{Classical finite element approximation}
	
	The structure is discretized with two-node elements using classical cubic Hermite shape functions \cite{chandrupatla2002} defined on the interval $\xi \in [-1,1]$:
	\begin{equation}\label{hermite_1d}
	W=\sum_{i=1}^2 \left[H^{\mathrm{w}}_i(\xi) \;\; -H^{\theta}_i(\xi)  h_x/2 \right]  \textbf{w}_i , 
	\end{equation}
	where 
    \begin{equation}\label{H_w}
	H^{\mathrm{w}}_i(\xi) =  (2+3\xi_i\xi-\xi_i\xi^3)/4,
	\end{equation}
	and
	\begin{equation}\label{H_theta}
    H^{\theta}_i(\xi) = (-\xi_i-\xi+\xi_i\xi^2+\xi^3)/4.
	\end{equation}
    Here, $\xi_{i=1,2}=\mp 1$ and $ \textbf{w}_i ^{\text{T}} = [ W_i \; \theta_{y,i} ]$ are the displacement and rotation at node $i$, respectively, with $\theta_y=-\partial W/ \partial x$. 
    Each element of length $h_x = x_2 -x_1$ is defined by the linear geometric mapping
    \begin{equation}\label{1d_mapping}
    x(\xi) = \sum_{i=1}^2 N_i x_i ,
    \end{equation}
    where $x_i$ are the nodal positions and $N_i = (1+\xi_i \xi)/2$ are the classical linear shape functions.

	\begin{figure}[htb!]
		\begin{center}
			\includegraphics[width=0.75\textwidth]{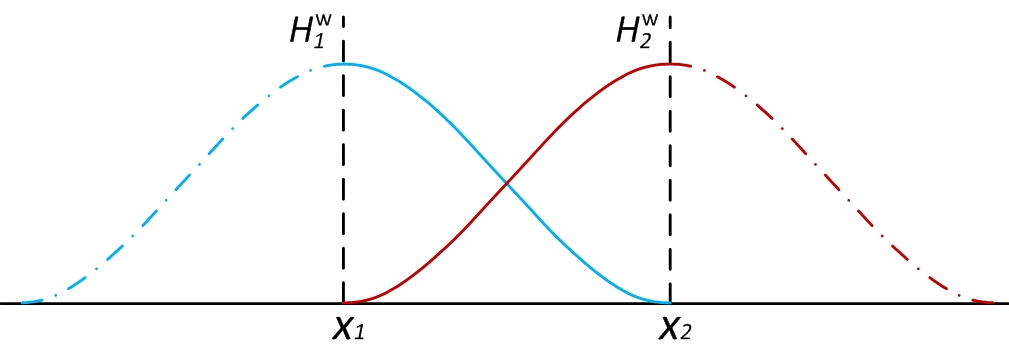}
		\end{center}
		\caption{Hermite-type PU functions}
		\label{fig_1d_show}
	\end{figure}
	
	\subsection{PUFEM approximation}	
	
    Classical PUFEM approximation which uses standard linear Lagrangian shape functions cannot be employed in the weak formulation \cite{zhou2019, babuvska2004}.
	One way to remedy this is to use the components of the Hermite functions associated with nodal displacement, i.e. $H^{\mathrm{w}}_{1,2}$ which form a partition of unity and, by construction, have a vanishing derivative at the endpoints as illustrated in Fig.~\ref{fig_1d_show}.
	The lateral displacement is thus expanded as follows
	\begin{equation} \label{pufem_1d}
	W =\sum_{i=1}^2 H^{\mathrm{w}}_i(\xi)  \sum_{n=1}^{N_i} A_i^n  \: \Phi_i^n.
	\end{equation}
	Here, enrichment functions $\Phi_i^n$ are assumed to be sufficiently smooth and Eq.\eqref{pufem_1d} can thus be regarded as a smooth partition of unity approximation which guarantees that first derivatives are continuous.
	Enrichment functions are normally chosen among the set of functions satisfying the partial differential equation, as advocated in Refs. \cite{Melenk, bathe2012}, but other types of functions, usually of polynomial and trigonometric form \cite{zhou2019, duarte2000, cornaggia2020, SHANGHSU2016}, can also be used. The set of polynomial functions is defined as 
	\begin{equation}\label{polyn_enrich_1d}
	\mathcal{P}_p =  \{1, \tilde{x}_i, \tilde{x}_i^{\;2}, \tilde{x}_i^{\;3}, \ldots  , \tilde{x}_i^{\;p}\},
	\end{equation}
	where $\tilde{x}_i= x-x_i$. 
	For a pure polynomial enrichment, the constant and linear terms correspond to the translational and rotational motions at nodes, respectively.
	A second type enrichment contains the free wave solutions to the governing equation (\ref{EoM_1d}) \cite{fahy2007, guyader2006}. In this work we only select propagating waves and avoid the use of evanescent waves as these functions are likely to have a detrimental effect on the conditioning of the algebraic system as frequency increases.
    Finally, enrichment functions are chosen among  the following set of functions
    \begin{equation}\label{wave_enrich_1d}
	\Phi_i^n \in  \{ \exp ( j k \tilde{x}_i ),\exp ( -j k \tilde{x}_i )  \} \cup \mathcal{P}_p,
	\end{equation}
	where 
	\begin{equation}\label{wavenumber}
	k=(\rho H \omega ^2/D)^{1/4}
	\end{equation}
	is the flexural wavenumber. 

    \subsection{Numerical examples}

    A concentrated harmonic point loading with a unit amplitude is applied at $L_x/4$ from the left end, corresponding to a node of the mesh. Geometrical parameters and mechanical properties of the plate are tabulated in Table \ref{properties}.
    
    	\begin{table}[htb!]
		\centering
		\begin{tabular}{c c}
		    \hline
			Length and thickness & $L_x=0.5$ m, $H=0.002$ m\\
			Material properties  & $E = 210$ GPa, $\rho = 7800$ kg/m$^3$, $\nu = 0.3$ \\
			\hline
		\end{tabular}
		\caption{Parameters used in our calculations.}
		\label{properties}
	\end{table}
    

   \begin{figure}[htb!]
		\begin{center}
		\includegraphics[width=0.48\textwidth]{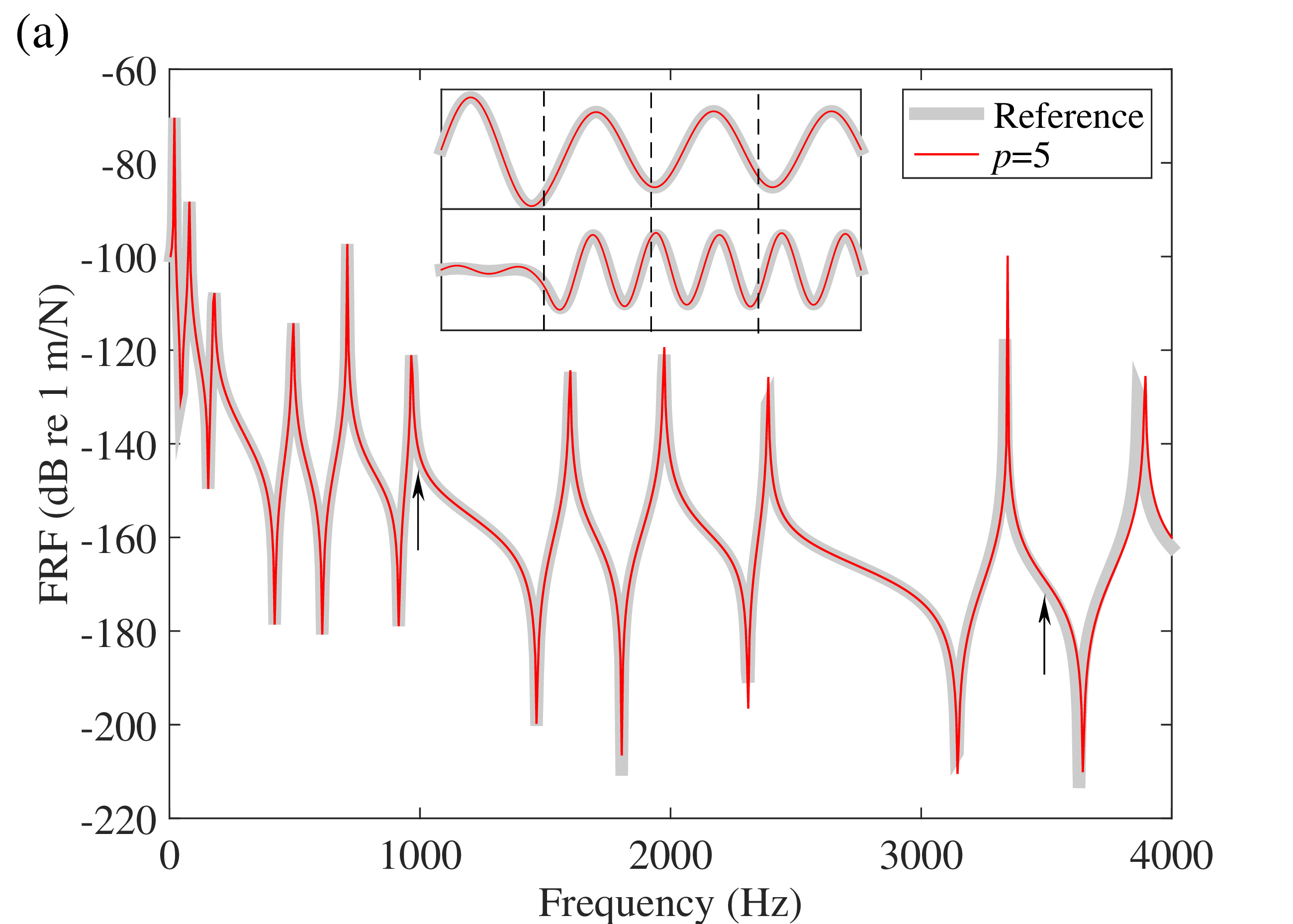}	\includegraphics[width=0.48\textwidth]{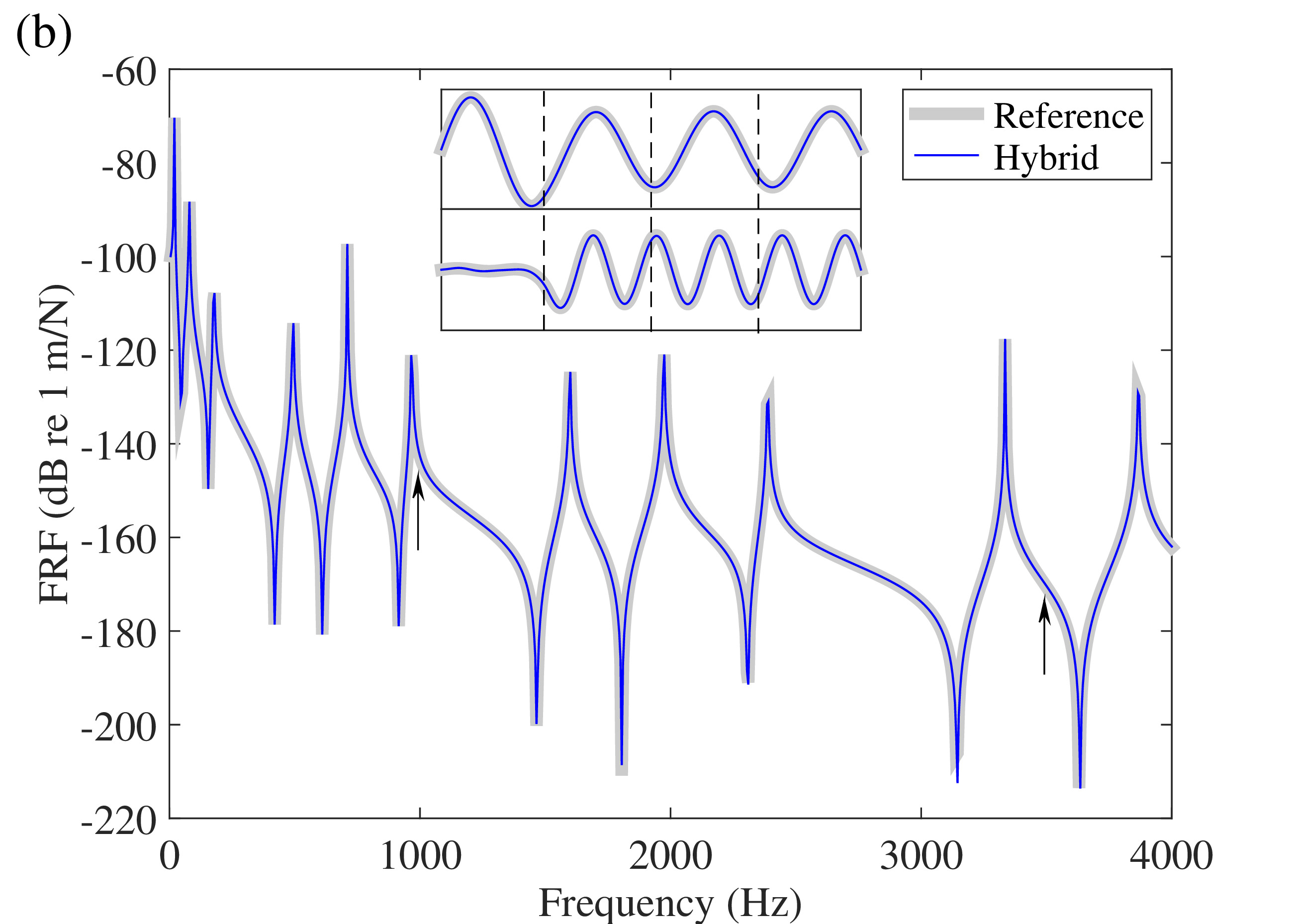}
		\end{center}
		\caption{Frequency response function (FRF) at the driving point of reference solutions and PUFEM with (a) quintic polynomial enrichment $p=5$ and (b) hybrid enrichment combining two progressive waves and cubic polynomials $p=3$ on 4 elements. The sub-figures are the corresponding structural responses at 1000Hz (upper) and 3500Hz (Bottom), located by arrows below FRFs, and the dash vertical lines are nodal positions.}
		\label{fig_1d_response}
	\end{figure}    
	
    Figure \ref{fig_1d_response} shows Frequency Response Function (FRF)  at the driving point using PUFEM with quintic polynomial enrichment $p=5$ and with a hybrid enrichment combining two waves propagating in opposite directions and a complete cubic polynomials $p=3$. 
    The PUFEM calculation is carried out with $M=4$ uniform elements of equal length. In this example, there are $N_i=6$ terms associated with each node of the PUFEM mesh so the two enrichment strategies are compared fairly. Here, a reference  solution is calculated using convergent modal series. 
    Clearly, there is a good agreement between the reference solution and the PUFEM results regardless of the type of enrichment used in the method.
    The corresponding structural responses along the $x$-axis are illustrated at 1000Hz and 3500Hz showing the capability of the method to simulate more than a wavelength over a single element.
    One can notice that the PUFEM quintic polynomials results show slight discrepancies above 3000 Hz where resonance peaks are shifted to a higher value, which is a typical phenomenon encountered when using conforming elements \cite{petyt2010}.
    The PUFEM with hybrid enrichment maintains a good accuracy even at higher frequencies, outperforming the pure polynomial enrichment with $p=5$.


    \begin{figure}[htb!]
		\begin{center}
		\includegraphics[width=0.48\textwidth]{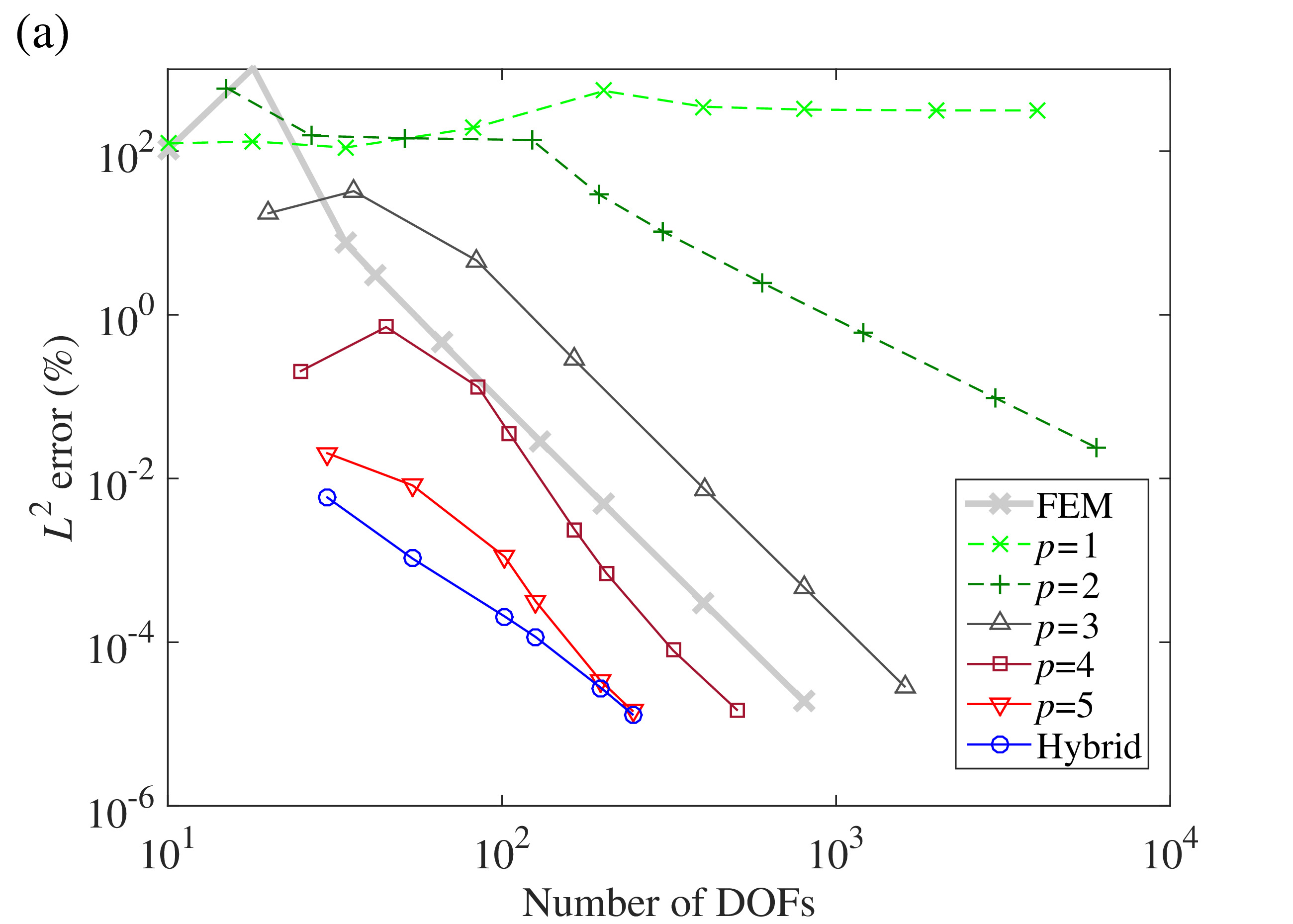}	\includegraphics[width=0.48\textwidth]{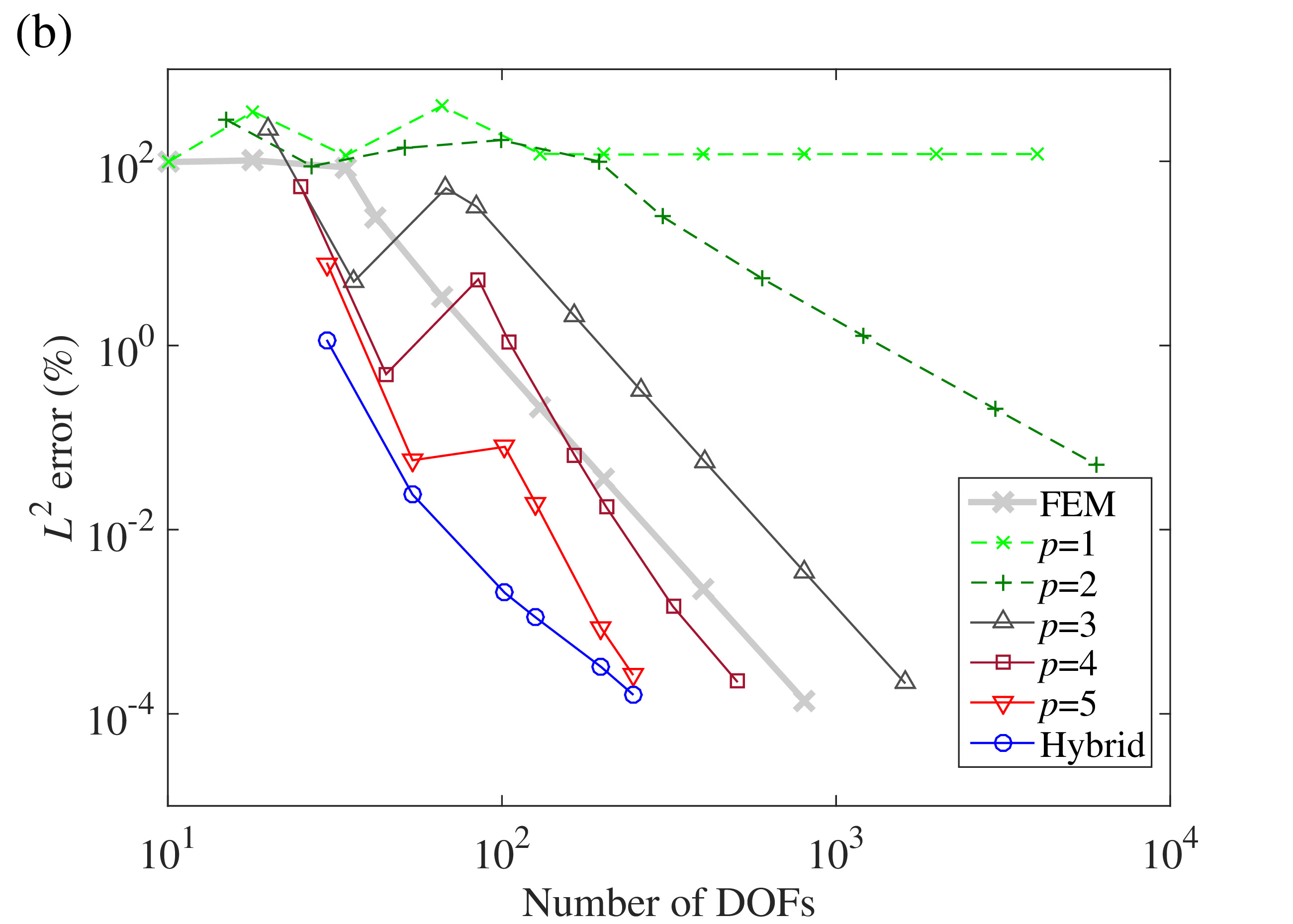}
		\end{center}
		\caption{Convergence curves with mesh refinement for classical FEM, the PUFEM polynomial enrichment with different orders $p$ and the hybrid PUFEM enrichment by combining the progressive waves and cubic polynomials $p=3$ at (a) 1000Hz and (b) 3500Hz.}
		\label{fig_1d_convergence}
	\end{figure}
	
	In order to illustrate this in a more systematic way, Fig.~\ref{fig_1d_convergence} shows a  convergence study of the classical FEM, the PUFEM with polynomial enrichments of order $p$, and the previous PUFEM with hybrid enrichment.
	These results  are obtained with a mesh refinement by increasing the number of elements of the same size starting with the coarse mesh with $M=4$.
	The relative $L^2$ errors are defined as
	\begin{equation}
	\varepsilon = \dfrac{\sqrt{ \int_0^{L_x} \vert W_{computed} - W_{ref} \vert^2 \mathrm{d}x} }{\sqrt{ \int_0^{L_x} \vert W_{ref} \vert^2 \mathrm{d}x } } \times 100 \% \; ,
	\end{equation}
	where $W_{ref}$ is the reference solution.  $L^2$ errors are plotted versus the total number of degrees of freedom. 
    Results show that the PUFEM with polynomial enrichment guarantee the convergence as long as $p\geq2$ which suggests that low order enrichment with constant and linear terms is not sufficient for the simulation of thin plate bending problem.
    Further analysis shows that errors behave like $\varepsilon \sim C h^{\sigma}$ where $\sigma \approx 2, 4$ for $p=2, 3$, respectively (here $C$ is a constant and the element length $h=h_x$ is inversely proportional to the number of DOFs). It can be observed that  classical FEM has the same convergence rate as the PUFEM with cubic polynomial enrichment and $\varepsilon \sim C' h^4$. The above analysis agrees well with theoretical estimates associated with high order finite element approximation for thin plate bending $O(h^{p+1}+h^{2(p-1)})$ \cite{strang1973, li2014}.
    Note that 8 polynomial terms (\ref{polyn_enrich_1d}) are associated with one PUFEM element enriched with cubic polynomials $p=3$ whereas classical element using complete cubic Hermite polynomials has only four DOFs; and although they have the same asymptotic rate of convergence, 
    classical FEM is more efficient in terms of data reduction.
    The higher order polynomial enrichment with $p=4, 5$ exhibits better convergence behavior and better data reduction.
    Convergence curves for the PUFEM polynomials enrichment all show similar peculiar behavior with a sudden drop of performance before reaching an asymptotic regime and this phenomenon also appears for other frequencies. The reason for this definitely deserves more thorough analyses, and the authors are reluctant to speculate in this paper as to its cause. 
    The incorporation of wave functions in the PUFEM formulation can further reduce the computational errors whilst avoiding the drop of performance.
    In the present case, the hybrid enrichment combining propagating waves and polynomials was found to deliver the best results and other numerical tests using evanescent waves did not show any improvements. Results at 3500 Hz show that convergence rate is somewhat not optimal as the number of DOFs increases and this stems from the nature of the solution in the vicinity of the load point. This point will be discussed in more details in the next section on two-dimensional case.
    
    

    

    \section{Two-dimensional solution}
    
	The normal displacement $W$ of the plate mid-surface in the $z$ direction obeys the equation of motion:
	\begin{equation}
	D \left( \frac{\partial ^4 W}{\partial x^4} + 2\frac{\partial ^4 W}{\partial x^2 \partial y^2} + \frac{\partial ^4 W}{\partial y^4} \right) - \omega^2 \rho H W = f_z
	\end{equation}
	where $f_z$ is a distributed transverse loading. A point force at $(x_F, y_F)$ can be expressed by using two-dimensional Dirac delta function as $f_z=F \delta (x-x_F,y-y_F)$.
	The plate governing equation can be written in an alternative compact matrix form \cite{zienkiewicz2013}
	\begin{equation}\label{eom_2dm}
	\bm{\mathcal{L}}^\text{T} \textbf{D} \bm{\mathcal{L}} W - \omega^2 \rho H W = f_z,
	\end{equation}
	where
	\begin{equation}
	\bm{\mathcal{L}}=\left[ \frac{\partial ^2 \;}{\partial x^2} \;\; \frac{\partial ^2 \;}{\partial y^2} \;\; 2\frac{\partial ^2 \;}{\partial xy}\right] ^\text{T} \; \text{and} \;\; \textbf{D}=D
	\begin{bmatrix}
	1 & \nu & 0 \\
	\nu & 1 & 0 \\
	0 & 0 & \; (1-\nu)/2
	\end{bmatrix}.
	\end{equation}
	Following the derivation of the previous section, the two-dimensional plate bending problem is formulated into a weak variational form and the corresponding essential boundary conditions are prescribed by introducing Lagrange multipliers.
	For simply supported plates studied in this work, i.e. with zero displacement specified at the boundary $\Gamma$, the associated variational formulation writes
	\begin{equation}\label{Variation_2d}
	\int_{\Omega} \left[ (\bm{\mathcal{L}} \delta \hspace{-0.2mm} W) ^\text{T} \textbf{D} \bm{\mathcal{L}} W - \omega^2 \delta \hspace{-0.2mm} W \rho h W   -\delta \hspace{-0.2mm} W f_z \right] \mathrm{d}x\mathrm{d}y - \int  _{\Gamma} \delta \hspace{-0.2mm} W \mathit{\Lambda} \mathrm{d}\Gamma  = 0,
	\end{equation}
	where $\Omega$ is the area of the plate mid-surface and $\mathit{\Lambda}$ is identified as the shear force along the border $\Gamma$. 
	The best way to deal with the boundary terms is to weakly enforce the essential conditions as \cite{Laghrouche2, Chazot2013_pufem1, Chazot2014_pufem2}: 
	\begin{eqnarray}\label{BC_weakform}
	\int  _{\Gamma} \delta \hspace{-0.2mm} \mathit{\Lambda}  W \mathrm{d}\Gamma = 0, \hspace{0.6cm} \forall (\delta \hspace{-0.2mm} \mathit{\Lambda}) .
	\end{eqnarray}
	This approach possesses the advantage of preserving the symmetry of the linear system and it has already been used in a PUFEM context to ease the treatment of the coupling conditions between two domains with distinct mechanical properties \cite{ Chazot2013_pufem1, Chazot2014_pufem2}.

	\begin{figure}[htb!]
		\begin{center}
			\includegraphics[width=0.75\textwidth]{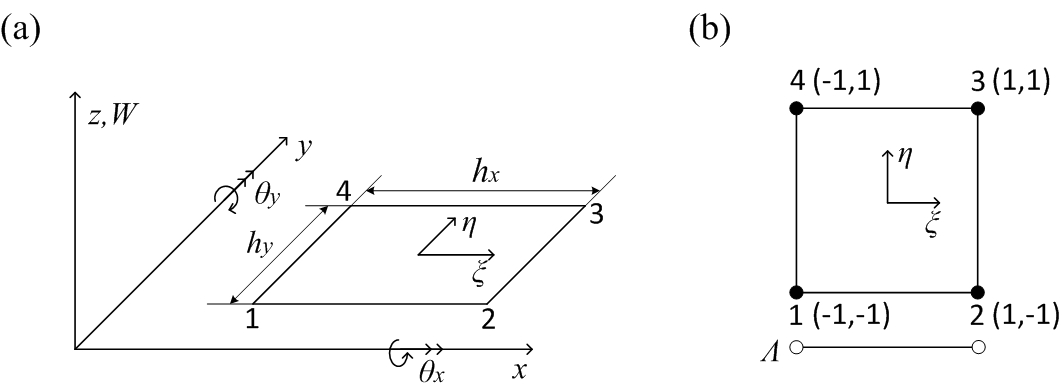}
		\end{center}
		\caption{A rectangular thin plate element and an edge element}
		\label{fig_2d_show}
	\end{figure}


	\subsection{Classical rectangular plate element (CR element)}
	
	A number of approaches are available in classical FEM to analyze thin plate flexural vibrations.
	One of the most common and effective ways is to use a two-dimensional conforming rectangular plate element called CR element \cite{bogner1965, petyt2010}.
	This element has four node points, with one situated at each corner, as shown in Fig.~\ref{fig_2d_show}. 
	Since the second derivatives for $W$ appear in the weak form (\ref{Variation_2d}), the DOFs and the shape functions attached to each node of the CR elements are devised to ensure that the displacement and its first derivatives with respect to $x$ and $y$ are continuous between elements.
	The nodal DOFs consist of the lateral displacement $W$, the rotations in two directions $\theta_y=-\partial W/ \partial x$ and $\theta_x= \partial W/ \partial y$, and the twist $W_{xy}= \partial^2 W/ \partial x \partial y$.
	A conforming CR element of dimensions $h_x \times h_y$ is constructed by taking the product of the unidimensional Hermite shape functions in Eq.(\ref{hermite_1d})
	\begin{equation}\label{cr_basis}
	W=\sum_{i=1}^4 \hat{\textbf{H}}_i(\xi,\eta)  \hat{\textbf{w}}_i 
	\end{equation}
	where
	\begin{equation}\label{cr_shape_fun}
	\hat{\textbf{H}}_i(\xi,\eta)=\left[ H^{\mathrm{w}}_i(\xi)H^{\mathrm{w}}_i(\eta) \hspace{0.2 cm} -H^{\theta}_i(\xi) H^{\mathrm{w}}_i(\eta)h_x/2  \hspace{0.2 cm} H^{\mathrm{w}}_i(\xi) H^{\theta}_i(\eta)h_y/2  \hspace{0.2 cm} H^{\theta}_i(\xi) H^{\theta}_i(\eta)h_x h_y/4 \right],
	\end{equation}
	and
	\begin{equation}
    \hat{\textbf{w}}_i ^{\text{T}} = [ W_i \hspace{0.4 cm}   \theta_{y,i}    \hspace{0.4 cm} \theta_{x,i}   \hspace{0.4 cm} W_{x y , i}   ] .
	\end{equation}
    The expressions for $H^{\mathrm{w}}_i$ and $H^{\theta}_i$ are given explicitly in Eqs. (\ref{H_w}) and (\ref{H_theta}). Here, $\xi,\eta \in [-1,1]$ are the local coordinates and $(\xi_i,\eta_i)$ are the local coordinates of nodes, see  Fig.~\ref{fig_2d_show}(b).
    CR elements are defined via the geometric mapping 
    \begin{equation}
	x = \sum_{i=1}^4 N_i (\xi,\eta) x_i \quad \text{and} \quad y = \sum_{i=1}^4 N_i (\xi,\eta) y_i,
	\end{equation}
	where $(x_i,y_i)$ are the location of the nodes of rectangular elements and $N_i= (1+\xi_i\xi)(1+\eta_i\eta)/4$.

	\subsection{PUFEM $C^1$ rectangular plate element}
	
\begin{figure}[htb!]
		\begin{center}
			\includegraphics[width=0.75\textwidth]{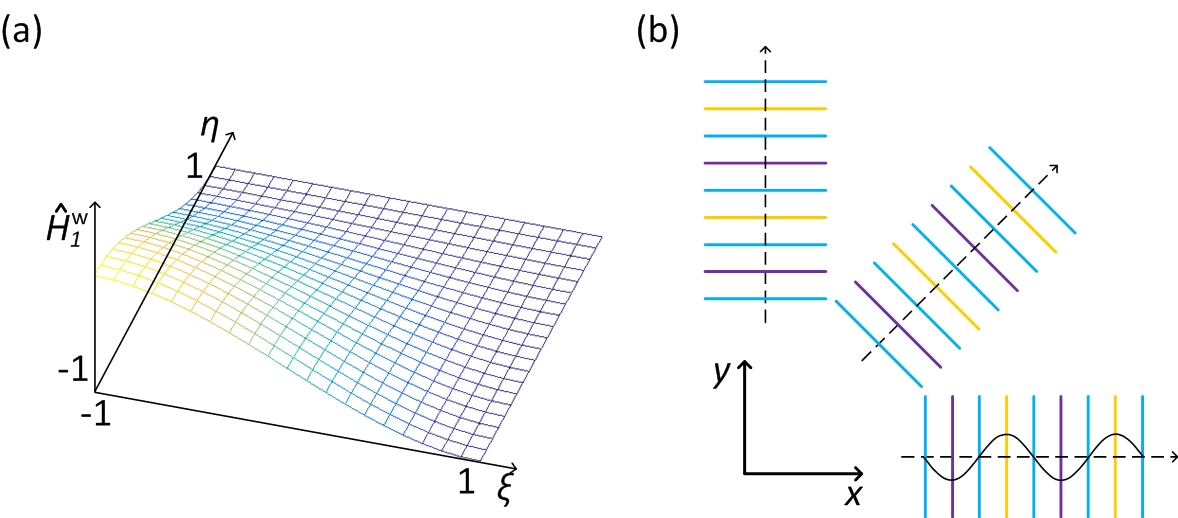}
		\end{center}
		\caption{(a) Hermite-type PU function $\hat{H}_1^{\mathrm{w}}$ and (b) plane flexural wave distributions for $\alpha_n=0, \pi/4, \pi/2$ over one quarter circle.}
		\label{fig_2d_wave}
	\end{figure}

	We can now mimic the developments given in the 1D formulation and select the classical CR element shape functions associated with nodal displacements in order to construct a partition of unity method with sufficient regularity.
	The plate transverse displacement field is expanded as follows
	\begin{equation}\label{pufem_2d_basis}
	W =\sum_{i=1}^4 \hat{H}_i^{\mathrm{w}}(\xi,\eta)  \sum_{n=1}^{\hat{N}_i} A_i^n  \: \Psi_i^n.
	\end{equation} 
	Here, functions $\hat{H}_i^{\mathrm{w}}(\xi,\eta)=H^{\mathrm{w}}_i(\xi)H^{\mathrm{w}}_i(\eta)$ correspond to the first term in (\ref{cr_shape_fun}), which is graphically illustrated in Fig.~\ref{fig_2d_wave}(a), and form a partition of the unity as $\sum_{i=1}^4\hat{H}_i^{\mathrm{w}}=1$,  $\forall (\xi,\eta)$ .
	By construction,  the approximation given by Eq. \eqref{pufem_2d_basis} guarantees that first derivatives with respect to both $x$ and $y$ directions are continuous everywhere and is therefore a good candidate for the discretization of the variational form.
	As for the enrichment strategies, we consider a polynomial enrichment using complete polynomials up to order $p$ in the two dimensions \cite{zienkiewicz2013, petyt2010}:
	\begin{equation}\label{polyn_enrich}
    \hat{\mathcal{P}}_p =  \{1, \tilde{x}_i, \tilde{y}_i, \tilde{x}_i^2, \tilde{x}_i\tilde{y}_i, \tilde{y}_i^2, \ldots  \},
	\end{equation}
	where $\tilde{x}_i= x-x_i$ and $\tilde{y}_i= y-y_i$. 
	The number of terms in the expansion is given by the Pascal triangle and $N_p=(p+1)(p+2)/2$.
	
	In light of the results obtained in the previous section, the hybrid enrichment tested in the present work consists in combining a set of plane waves propagating in different directions \cite{fahy2007, guyader2006} with a polynomial enrichment up to a given order $p$:
	\begin{equation}\label{wave_enrich}
	\Psi_i^n \in  \{  \ldots, \exp[j( k_{x,n}\tilde{x}_i +  k_{y,n}\tilde{y}_i) ] , \ldots\} \cup \hat{\mathcal{P}}_p,
	\end{equation}
	where
	\begin{equation}\label{wave_direction}
	[(k_{x,n})^2+(k_{y,n})^2]^2=\rho H \omega ^2/D=k^4.
	\end{equation}
	Here, the real-valued wavenumber components $k_{x,n}$ and $k_{y,n}$ correspond to a plane wave propagating with an angle $\alpha_n$ with respect to the $x$-axis and
	\begin{equation}\label{wave_angle}
	(k_{x,n}, k_{y,n})= k ( \cos{\alpha_n},  \sin{\alpha_n}), \; \text{where} \; \alpha_n=2\pi n /q_{i}, \; n = 1, \ldots, q_{i}.
	\end{equation}
	The flexural wave propagation directions $\alpha_n$, shown in Fig.~\ref{fig_2d_wave}(b), are chosen to be evenly distributed over the unit circle and $q_i$ corresponds to the number of plane waves attached to node $i$. Thus, the number of functions in the summation in \eqref{pufem_2d_basis} is $\hat{N}_i=q_i+N_p$.
	Note that the hybrid enrichment degenerates naturally to a pure polynomial enrichment by simply taking $q_i=0$.
	As opposed to the 1D case whereby only 2 plane waves exist, a refinement of the PUFEM element is achieved by increasing the number of plane waves attached locally to each node of the PUFEM mesh as in Ref. \cite{Chazot2013_pufem1, Chazot2014_pufem2} for the simulation of acoustic waves. Here again, the use of evanescent waves, see Ref. \cite{guyader2006} for more details, is not favored due to the presence of fast decaying components which can have a detrimental effect on the conditioning of the resultant system matrices. 
	Furthermore, it will be shown in Section 4 that the combination with low order polynomial functions is sufficient to capture the presence of evanescent waves which are excepted to only appear in the vicinity of a point load and near the edges of the plate.


	\subsection{Approximation for Lagrange multipliers}
	The Lagrange multiplier $\mathit{\Lambda}$, defined on the  border $\Gamma$ in (\ref{Variation_2d}), see Fig.~\ref{fig_2d_show}(b), is employed to prescribe constraints for both classical FEM and PUFEM in this work.
	The Lagrange multipliers over a two-nodes  edge element of the PUFEM mesh are approximated using polynomial-enriched elements
	similar to the 1D case. 
	For instance, the variable $\mathit{\Lambda}$ along a border aligned in the $x$-direction is expanded as
	\begin{equation} \label{lm_1d}
	\mathit{\Lambda} =\sum_{i=1}^2 H^{\mathrm{w}}_i(\xi)  \sum_{l=1}^{N_i} C_i^l  \: \Phi_i^l.
	\end{equation}
	where function $\Phi_i^l \in \mathcal{P}_p$ is chosen among the set of 1D polynomial enrichment functions \eqref{polyn_enrich_1d}.  The present discretization scheme for Lagrange multipliers is adopted here in order to ensure a sufficiently accurate description for the shear force on the boundary whilst avoiding ill-conditioning problems that might arise with more sophisticated functions.
	After substitution we arrive at a symmetric system of the 
	following form 
	\begin{equation}
	\begin{bmatrix}
	\textbf{K}_{\small{WW}} & \textbf{K}_{\small{W\mathit{\Lambda}}} \\
	\textbf{K}_{\small{W\mathit{\Lambda}}}^\text{T} & \textbf{0} \\
	\end{bmatrix}
	\begin{Bmatrix}
          \textbf{A} \\
          \textbf{C}
    \end{Bmatrix}
	= 
	\begin{Bmatrix}
          \textbf{F} \\
          \textbf{0}
    \end{Bmatrix}.
	\end{equation}
	where vectors $\textbf{A}$ and $\textbf{C}$ contain the unknown expansion coefficients in (\ref{pufem_2d_basis}) and  (\ref{lm_1d}) and $\textbf{F}$ is the loading vector.

	\section{Numerical results and discussions}
	
	The configuration for evaluating the performance of the PUFEM elements is a square flat plate of size $[0, L] \times [0, L]$ with all four edges simply supported. For any arbitrary load distribution $f_z$, this simple  configuration has a reference solution that can be built using the modal superposition method \cite{fahy2007, guyader2006}.
	The mechanical properties in Table \ref{properties} are used in the  calculations.
	In order to simplify the analysis, the same number of plane waves is attached to each node of the mesh and we put $q_i=q$. A small initial deflection angle $\pi/50$ is added to all propagation angles $\alpha_n$ in (\ref{wave_angle}), so that  wave propagation directions are not aligned with  the principle axes of the square plate. Similarly, for a given  order $p$, the number of terms $N_p$ in the polynomial expansion is chosen to be identical for all nodes.
	The number of polynomial terms in \eqref{lm_1d} is also identical for all nodes lying on the edges
	and is chosen to be in accordance with order $p$ so we take  $N_i=p+p'$, where $p'$ is an integer.

	\begin{figure}[htb!]
		\begin{center}
			\includegraphics[width=0.48\textwidth]{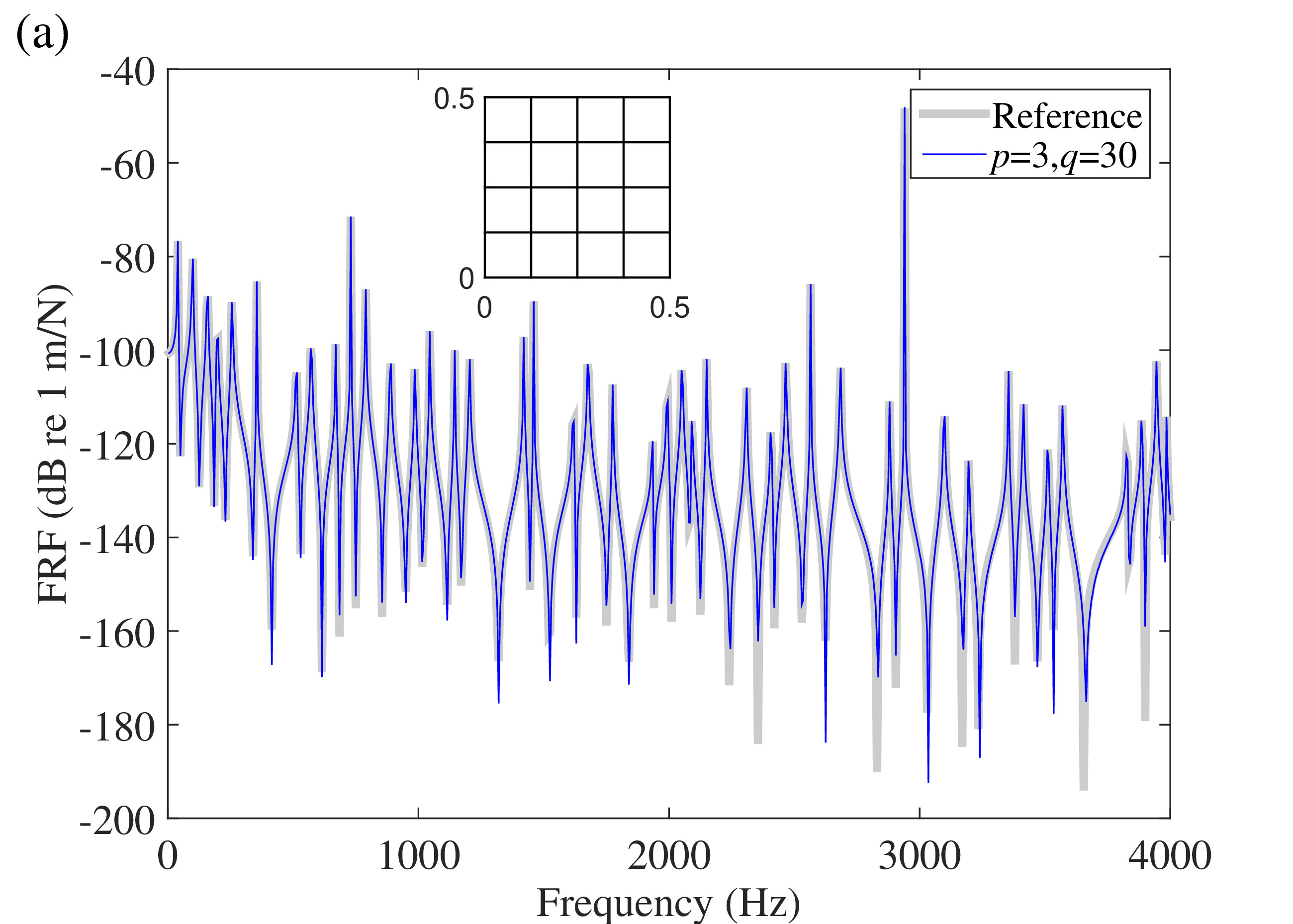}
			\includegraphics[width=0.48\textwidth]{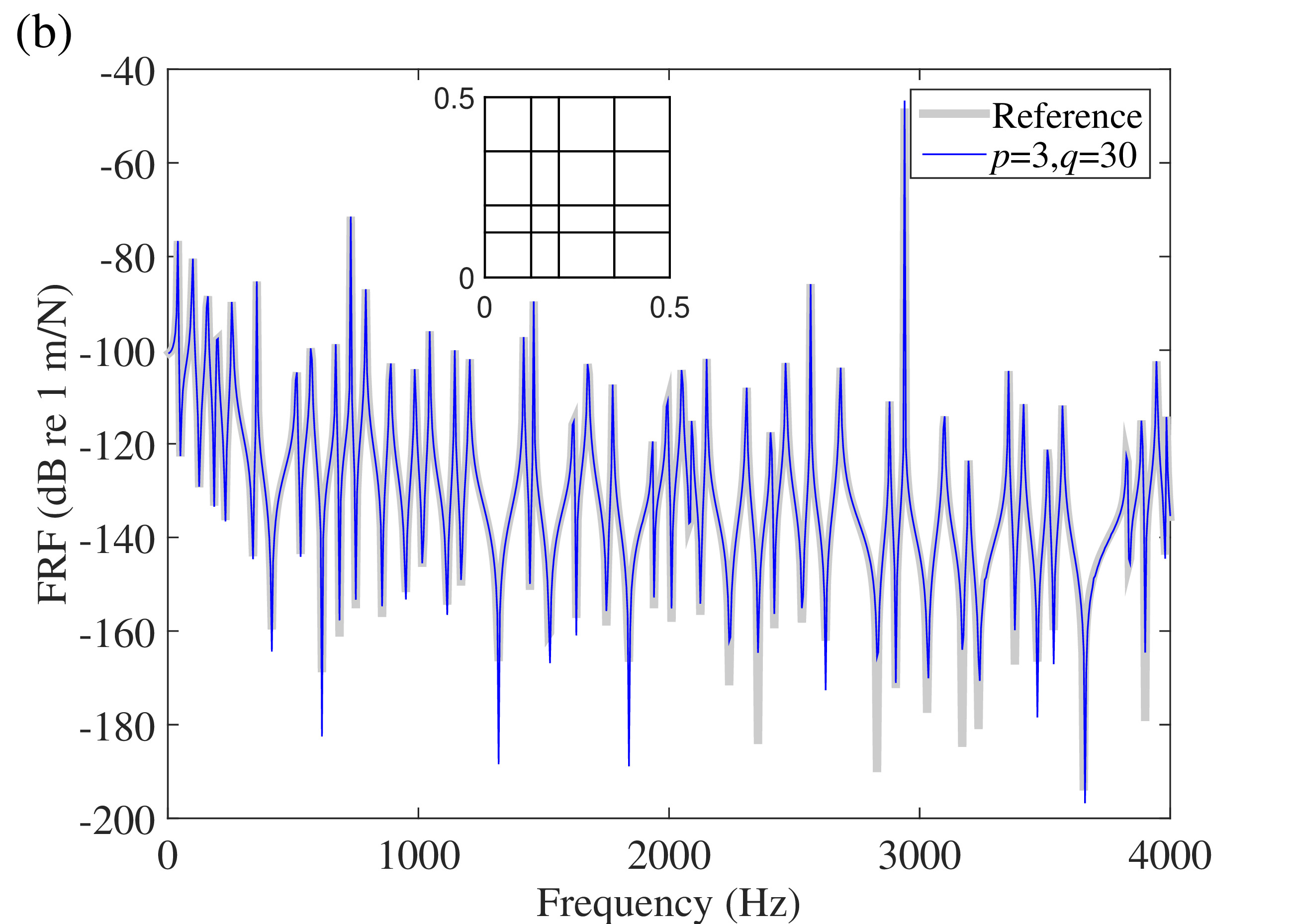}
			\includegraphics[width=0.48\textwidth]{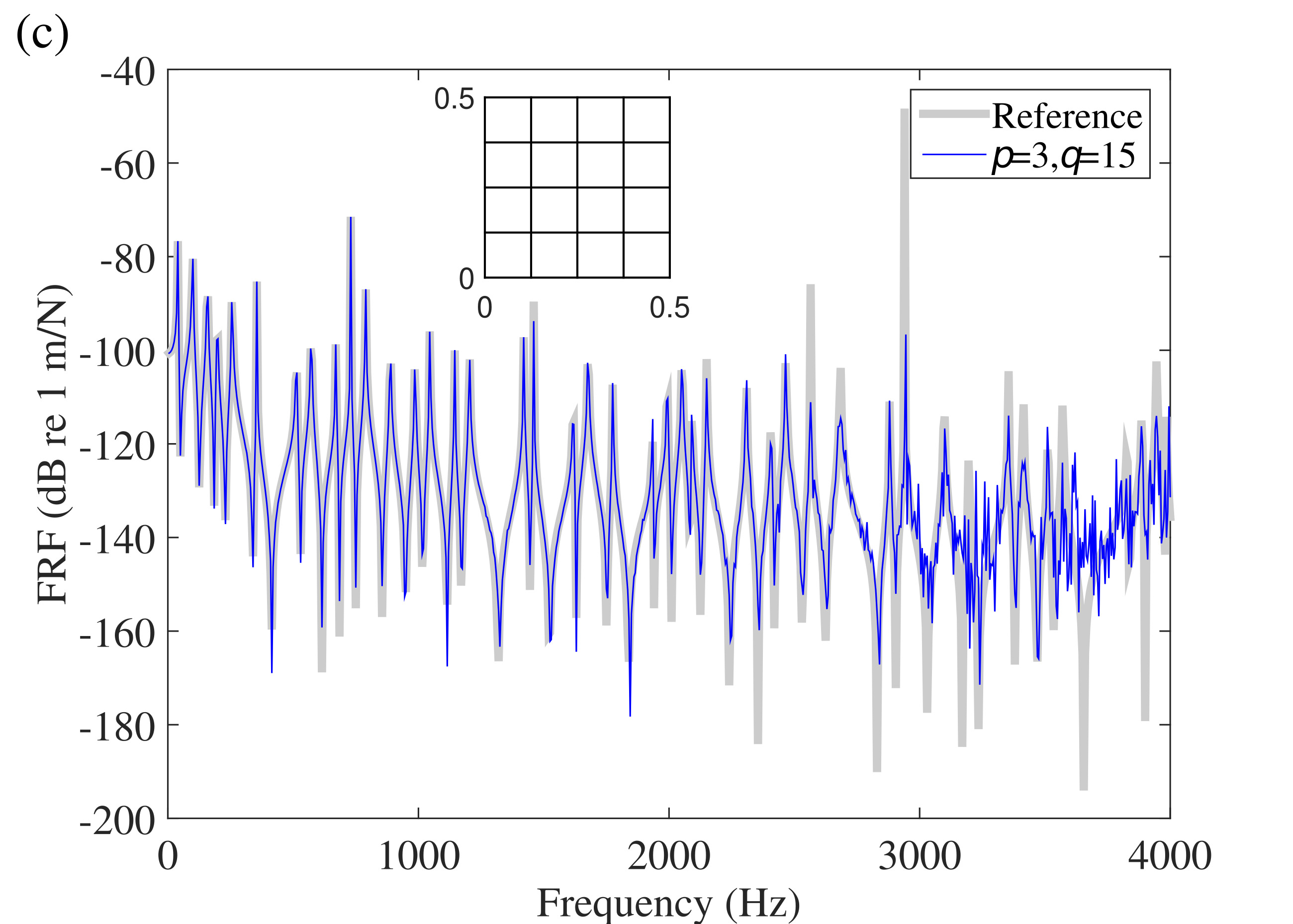}
			\includegraphics[width=0.48\textwidth]{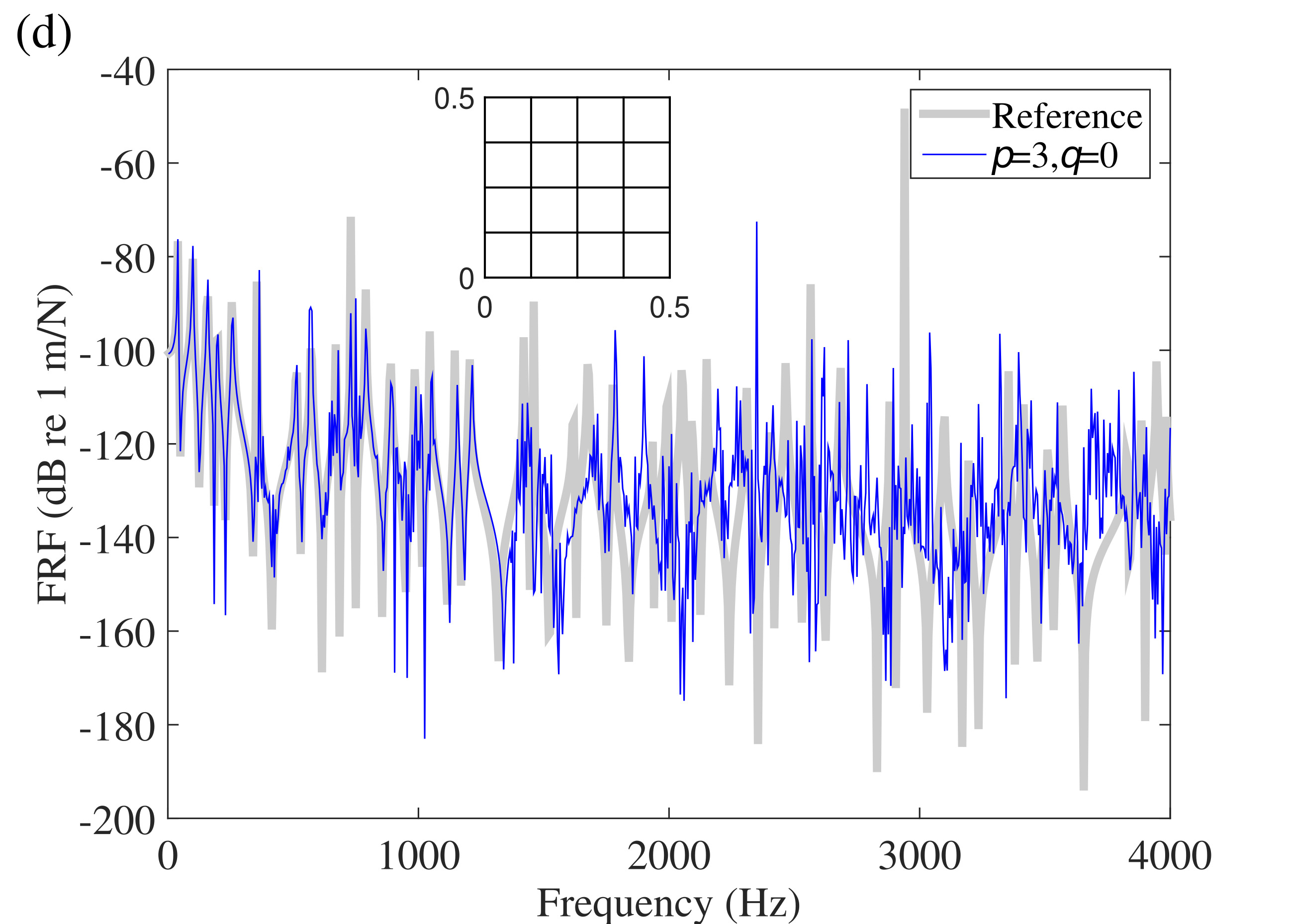}
		\end{center}
		\caption{Comparisons of driving point FRF between reference solutions and PUFEM with different numbers of free waves $q$ plus complete cubic polynomials $p=3$, using uniform mesh (a) (c) (d) and nonuniform mesh (b). Sub figures are the mesh grids in $x$-$y$ plane}
		\label{fig_2d_point}
	\end{figure}

	\subsection{Frequency response}
	
	The uniform mesh has 16 square elements of the same size. A nonuniform mesh is also investigated with PUFEM nodes defined on the grid $\{0, 0.25L, 0.4L, 0.7L, L\}^2$  (see Fig.~\ref{fig_2d_point} (b)).
    A concentrated harmonic point loading with a unit amplitude, $f_z=\delta(x-x_F, y-y_F)$, is applied at $(x_F,y_F)=(0.25L,0.25L)$, which corresponds to a node of the PUFEM mesh. The FRFs at the driving point are shown in Figs.~\ref{fig_2d_point} (a), (b), (c) and (d). 
	The structural responses are calculated with $q=0$, 15 and 30 plane waves combined with cubic polynomials, i.e. $p=3$.  From preliminary tests, numerical results with acceptable accuracy are obtained when $N_i$ is chosen above $p+2$ and numerical convergence was reached as soon as $N_i=p+8$, which is used in our calculations.
	For the special case of a pure polynomial enrichment ($q=0$), the polynomial order used for the  Lagrange multipliers is set to be the same as for the displacement, which means that $p'=0$.
	
\begin{figure}[htb!]
		\begin{center}
			\includegraphics[width=0.46\textwidth]{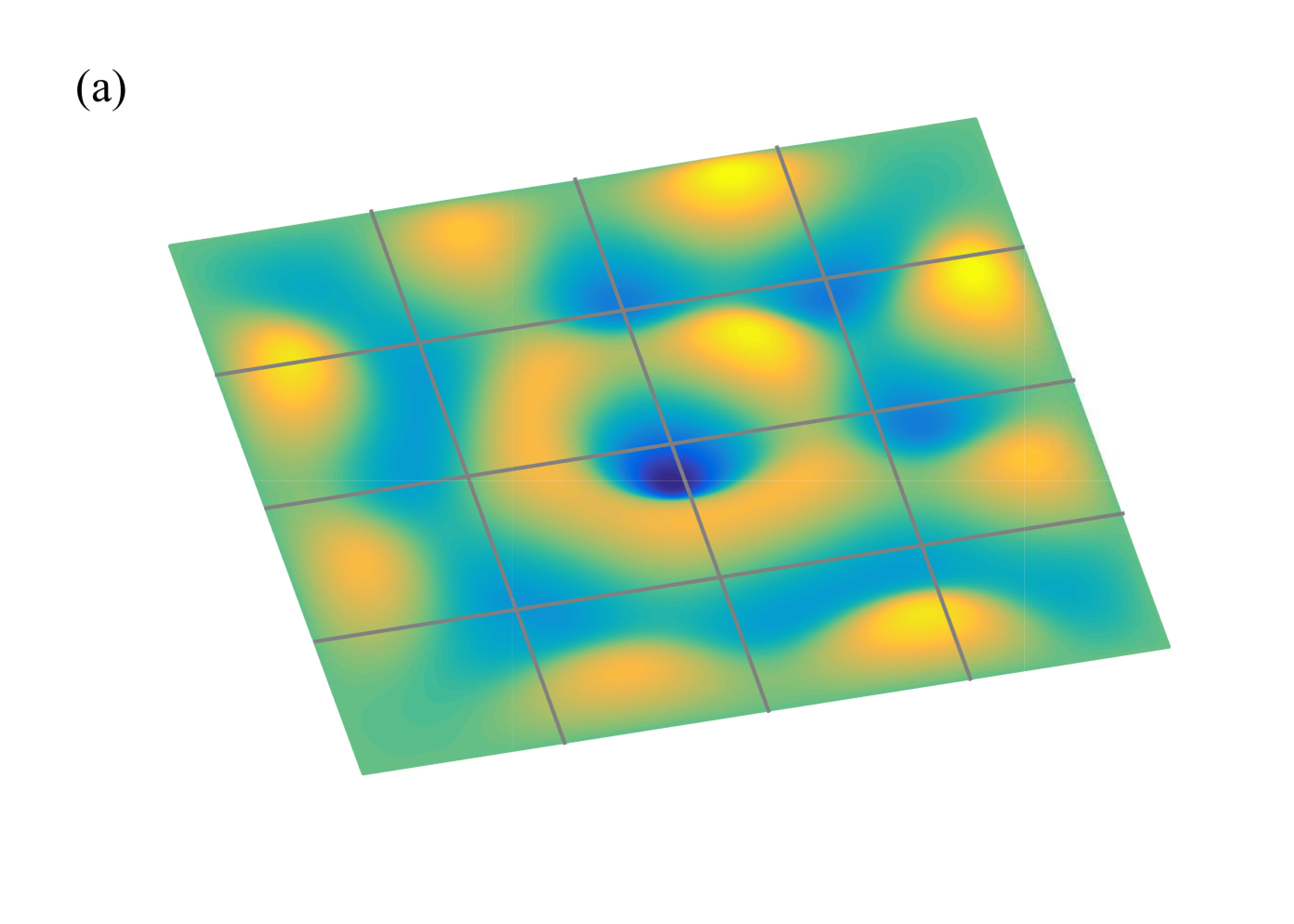}	\includegraphics[width=0.46\textwidth]{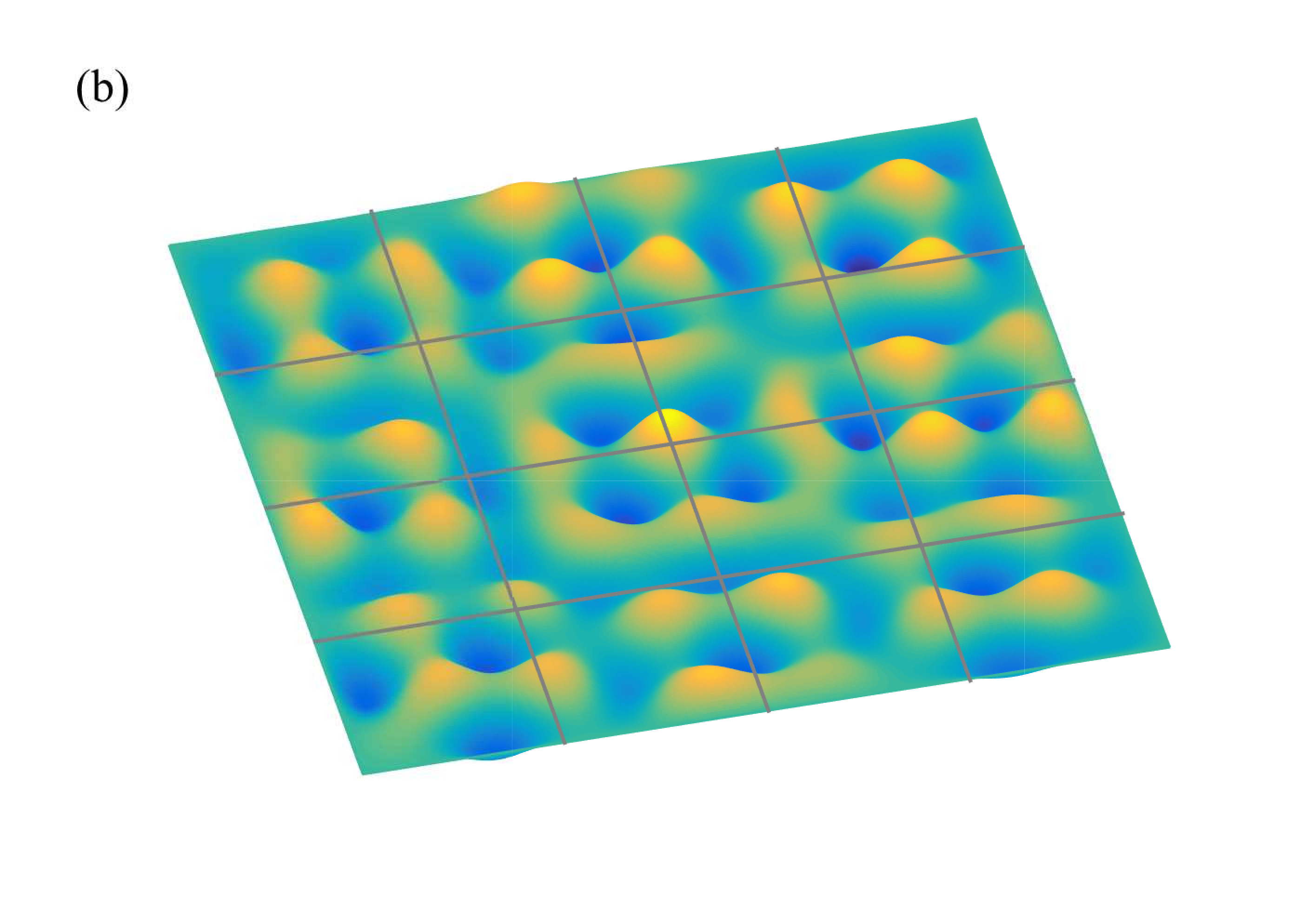}	\includegraphics[width=0.46\textwidth]{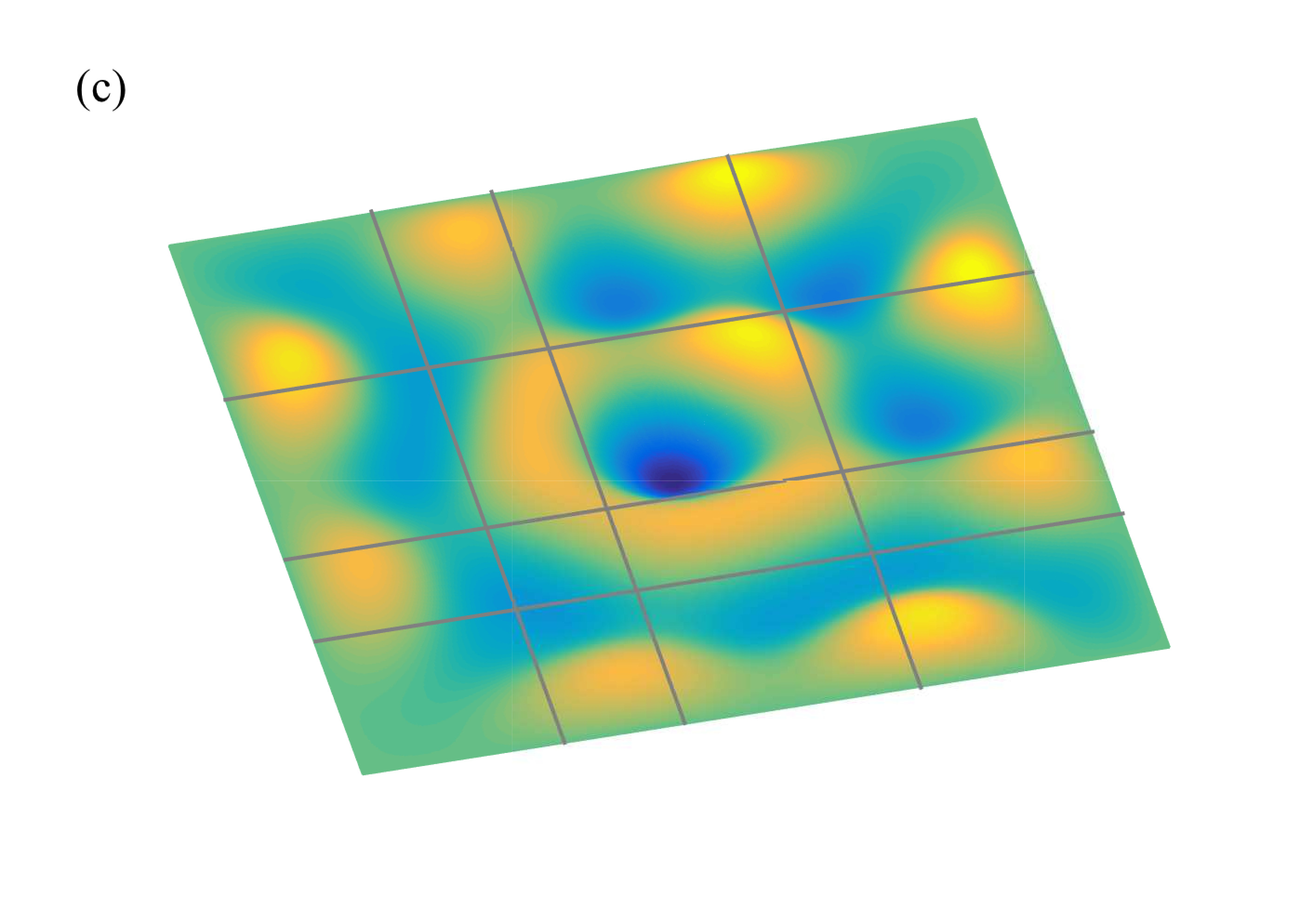}	\includegraphics[width=0.46\textwidth]{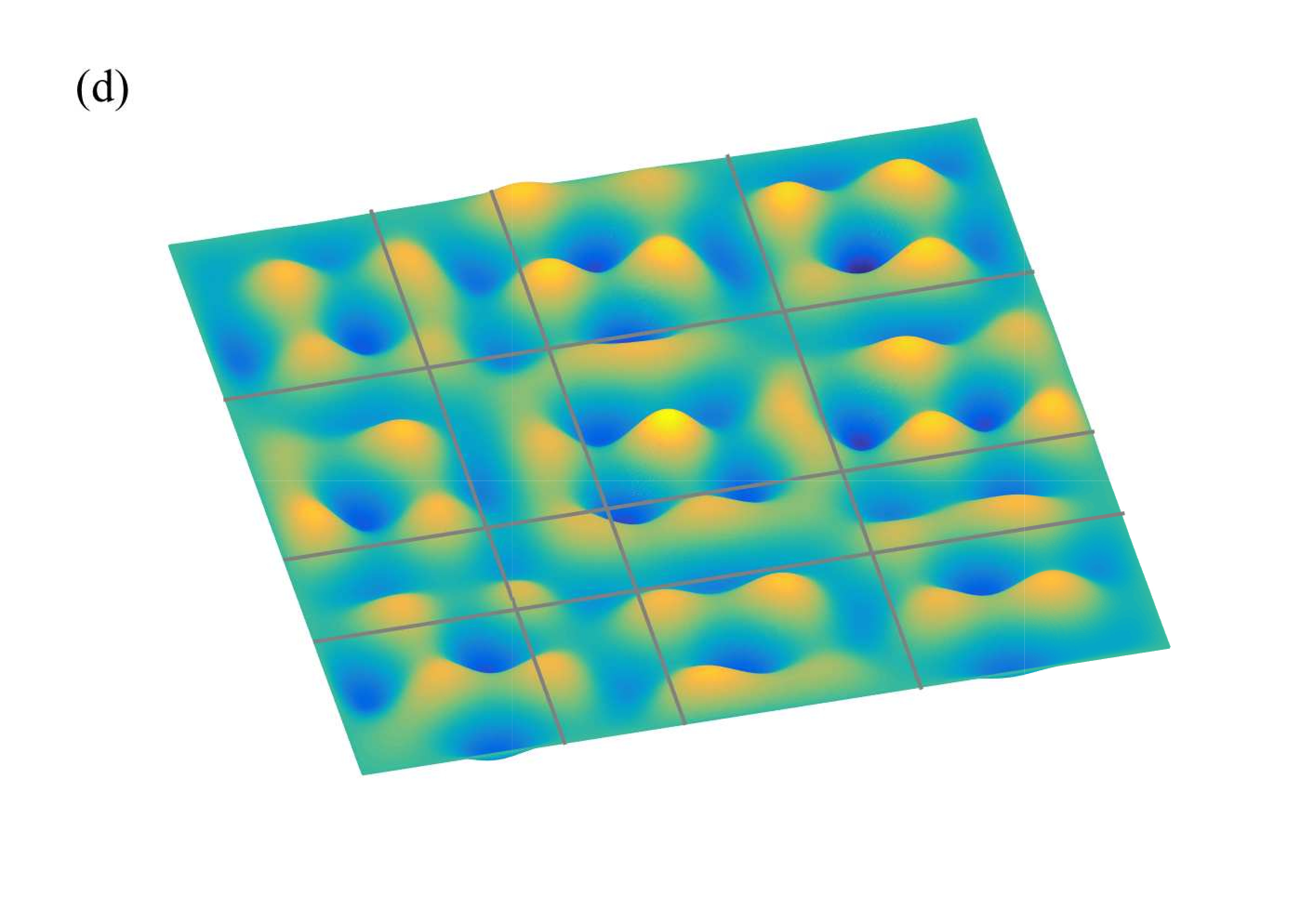}	\includegraphics[width=0.46\textwidth]{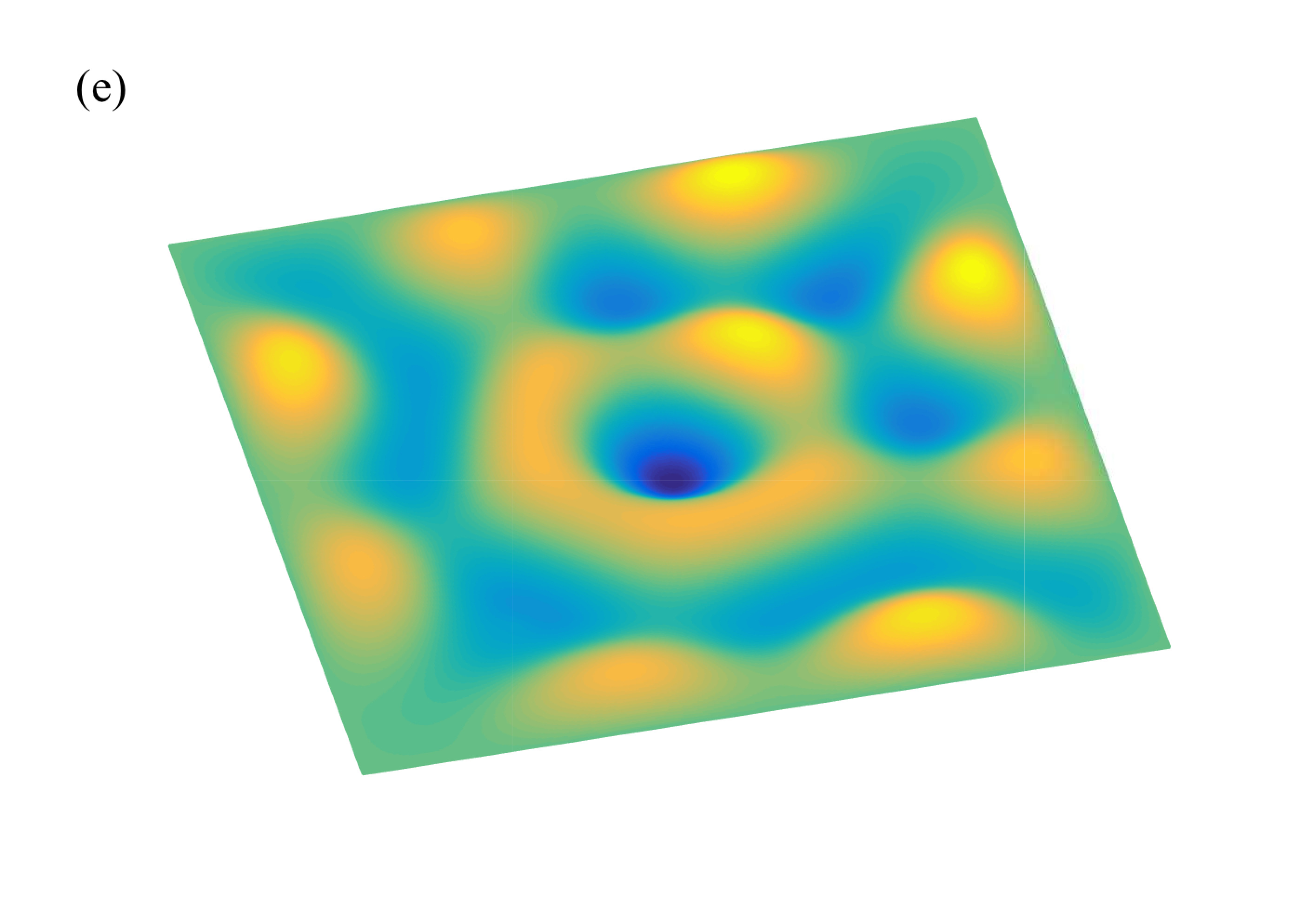}			\includegraphics[width=0.46\textwidth]{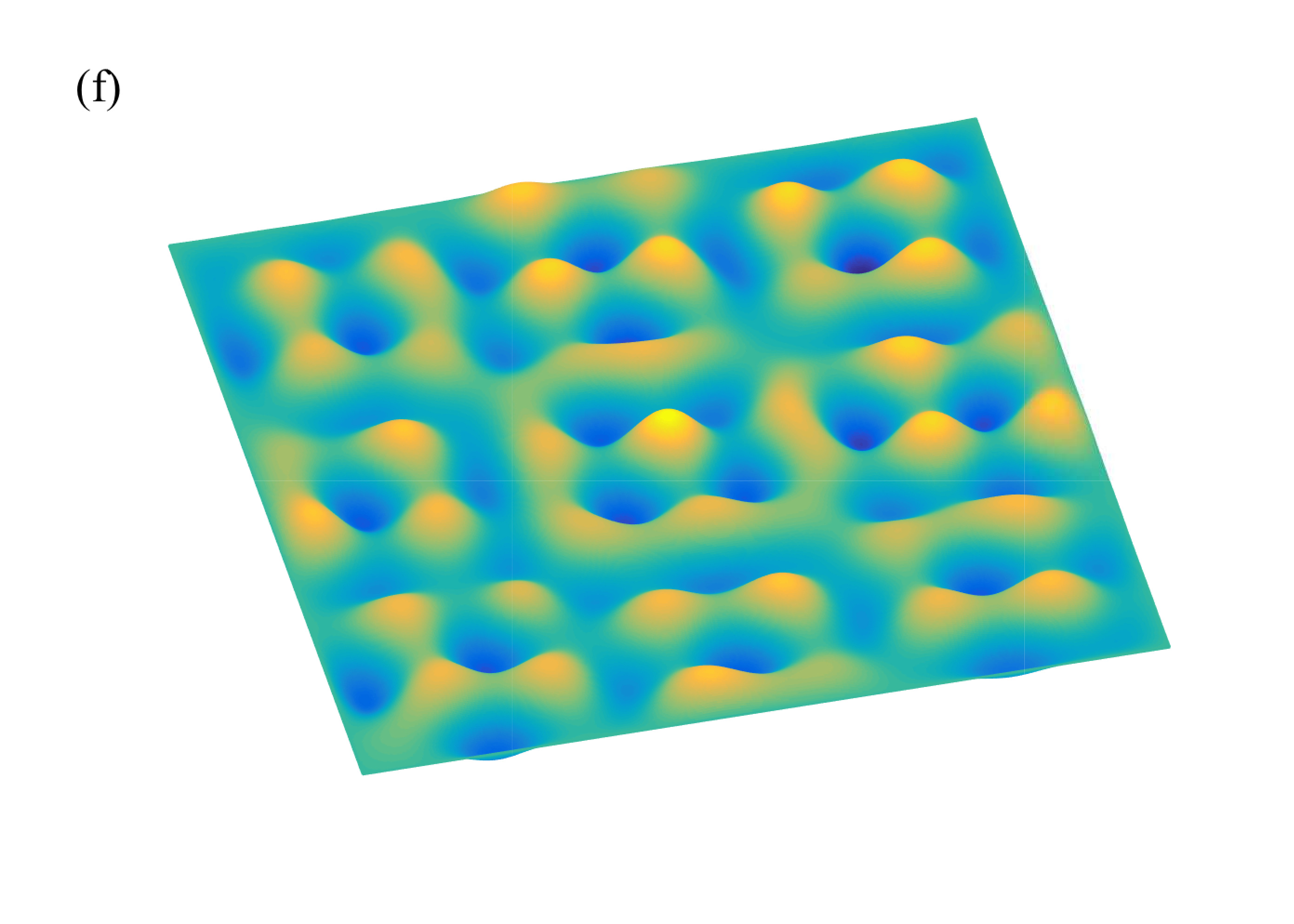}
			\includegraphics[width=0.4\textwidth]{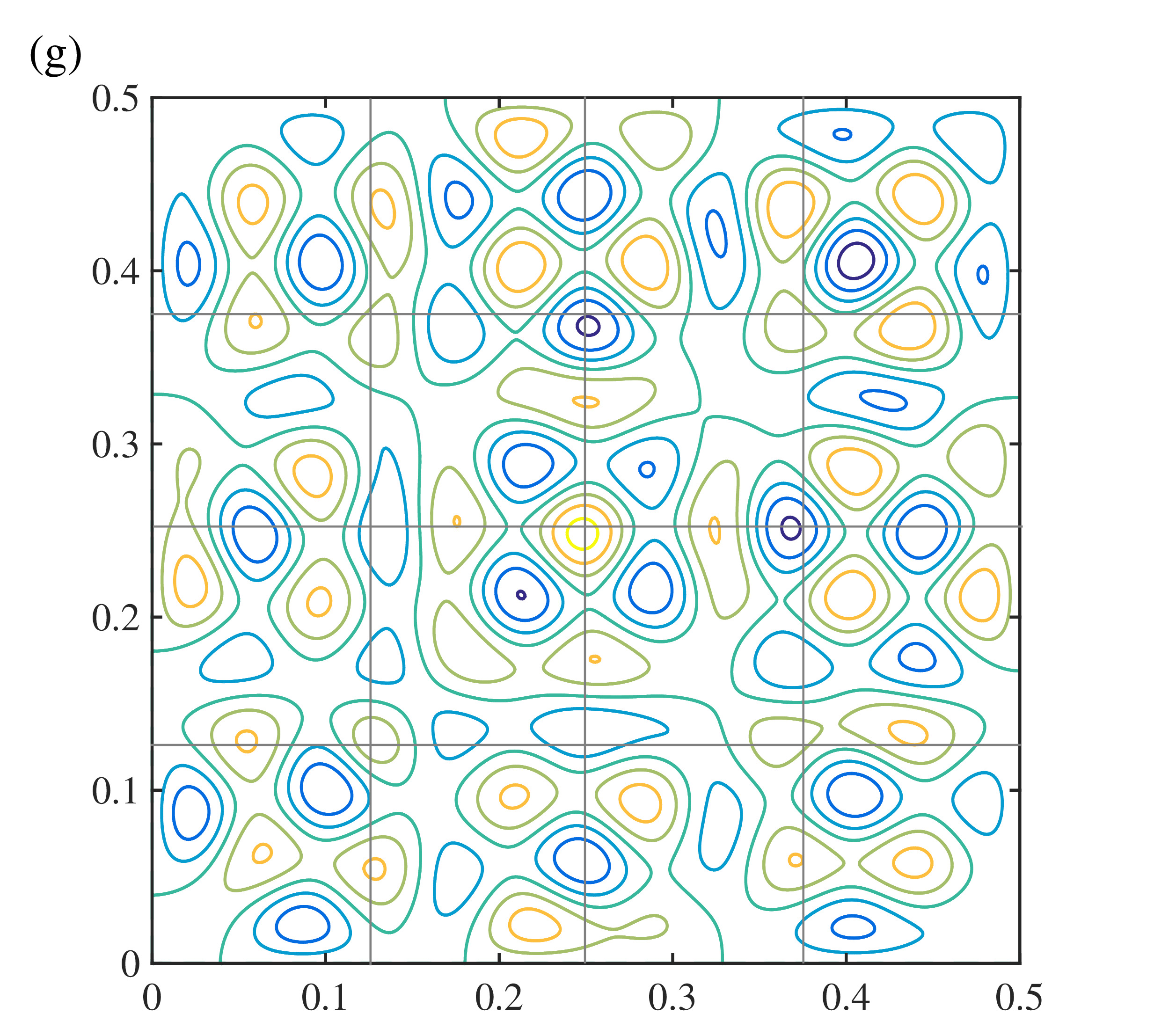}\hspace{0.2cm}
			\includegraphics[width=0.4\textwidth]{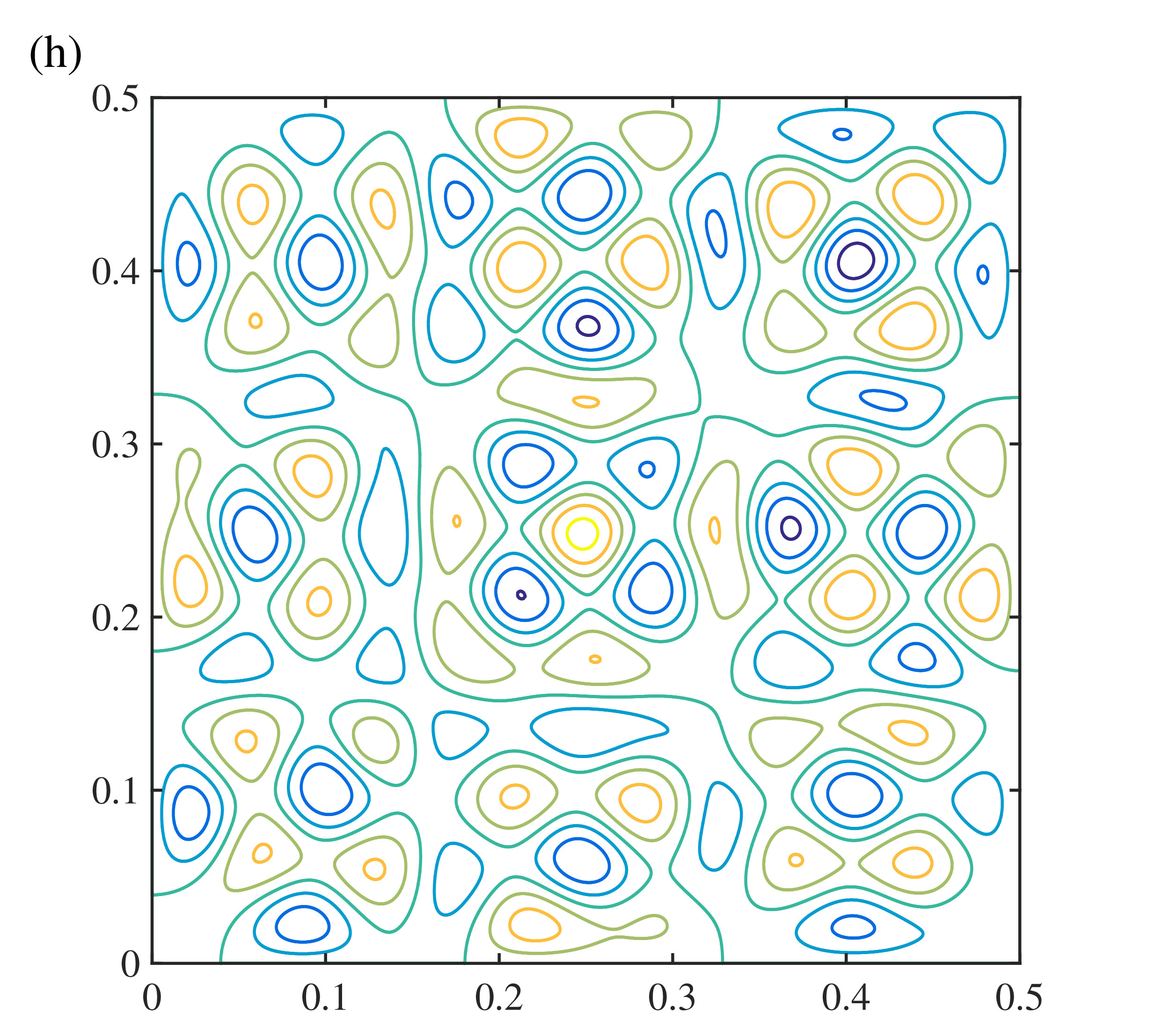}
		\end{center}
		\caption{Plate response calculated by PUFEM with 30 free waves plus cubic polynomials ($q=30, p=3$), using uniform mesh (a), (b) and nonuniform mesh (c), (d), and the corresponding reference solutions (e), (f). Vibration frequencies of (a), (c) and (e) are 1000Hz while others are 3500Hz. (g) and (h) are the contour plots of (b) and (f), respectively.}
		\label{fig_2d_response}
	\end{figure}

	Results show that the hybrid enrichment with $q=30$ plane waves delivers excellent results regardless of the PUFEM mesh. Differences are hardly noticeable. The reason for this lies in the fact that the enrichment basis is sufficient to simulate the wave field over large elements, here from $h=0.15L$ to $h=0.3L$ for the non-uniform mesh, at the upper frequency of interest, here $4000$Hz. 
	The corresponding plate displacement fields at two specific frequencies, 1000Hz and 3500Hz, are shown in Fig.~\ref{fig_2d_response}. The structural response computed by the PUFEM on both meshes agrees very well with the reference solutions which show symmetry with respect to the line $y=x$, as expected.
	As shown by results of Fig. \ref{fig_2d_point} (c), the wave-polynomial enrichment with $p=3, q=15$
	is not sufficient to capture the strong oscillating behaviors over the PUFEM element when the wavelength becomes too short, here around $2500$ Hz, and discrepancies appear at higher frequencies.
	The pure polynomial enrichment, see Fig. \ref{fig_2d_point} (d), is only viable at low frequencies and in this scenario a mesh refinement is necessary.

	\subsection{Convergence studies}

\begin{figure}[htb!]
		\begin{center}
			\includegraphics[width=0.48\textwidth]{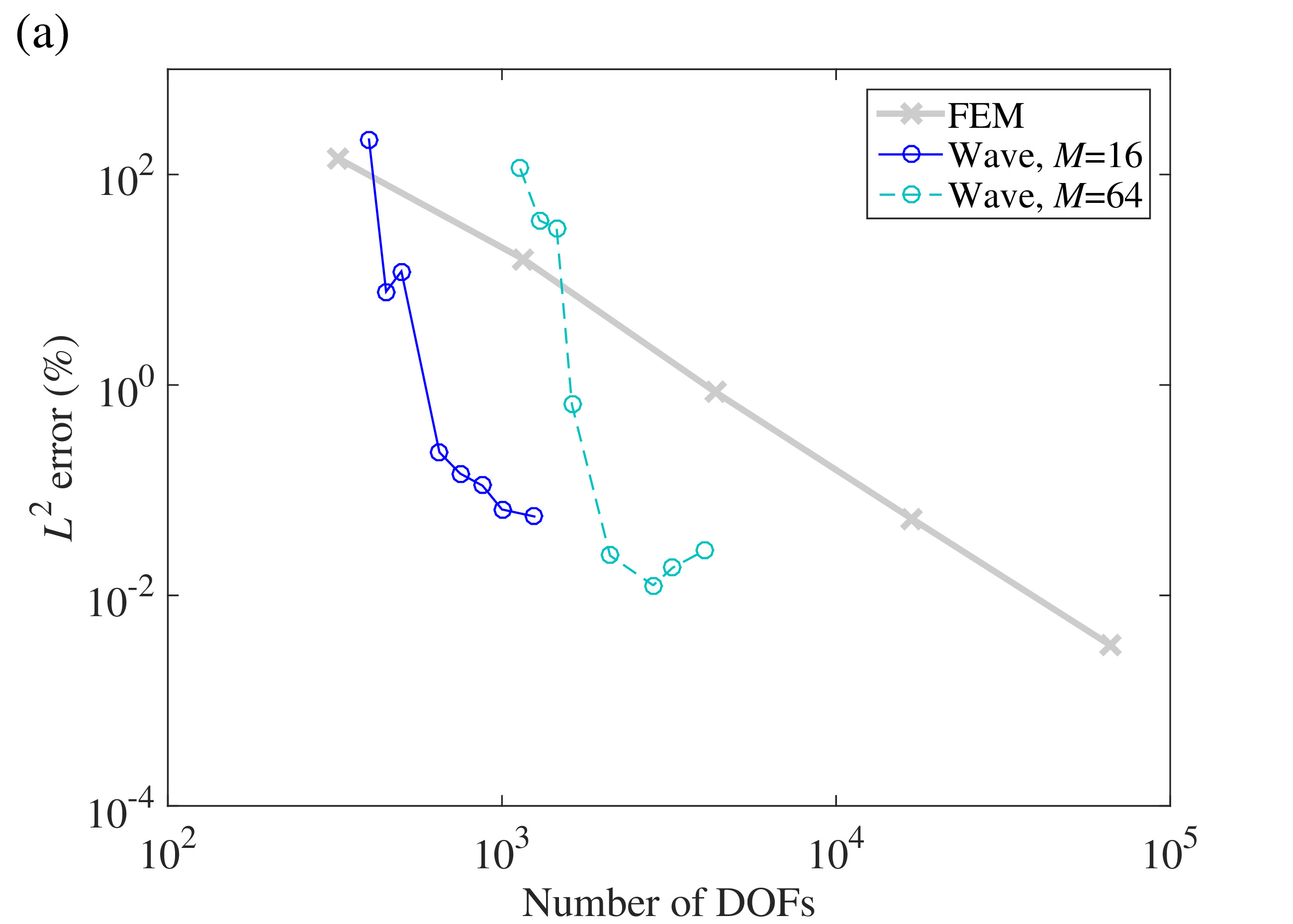}
			\includegraphics[width=0.48\textwidth]{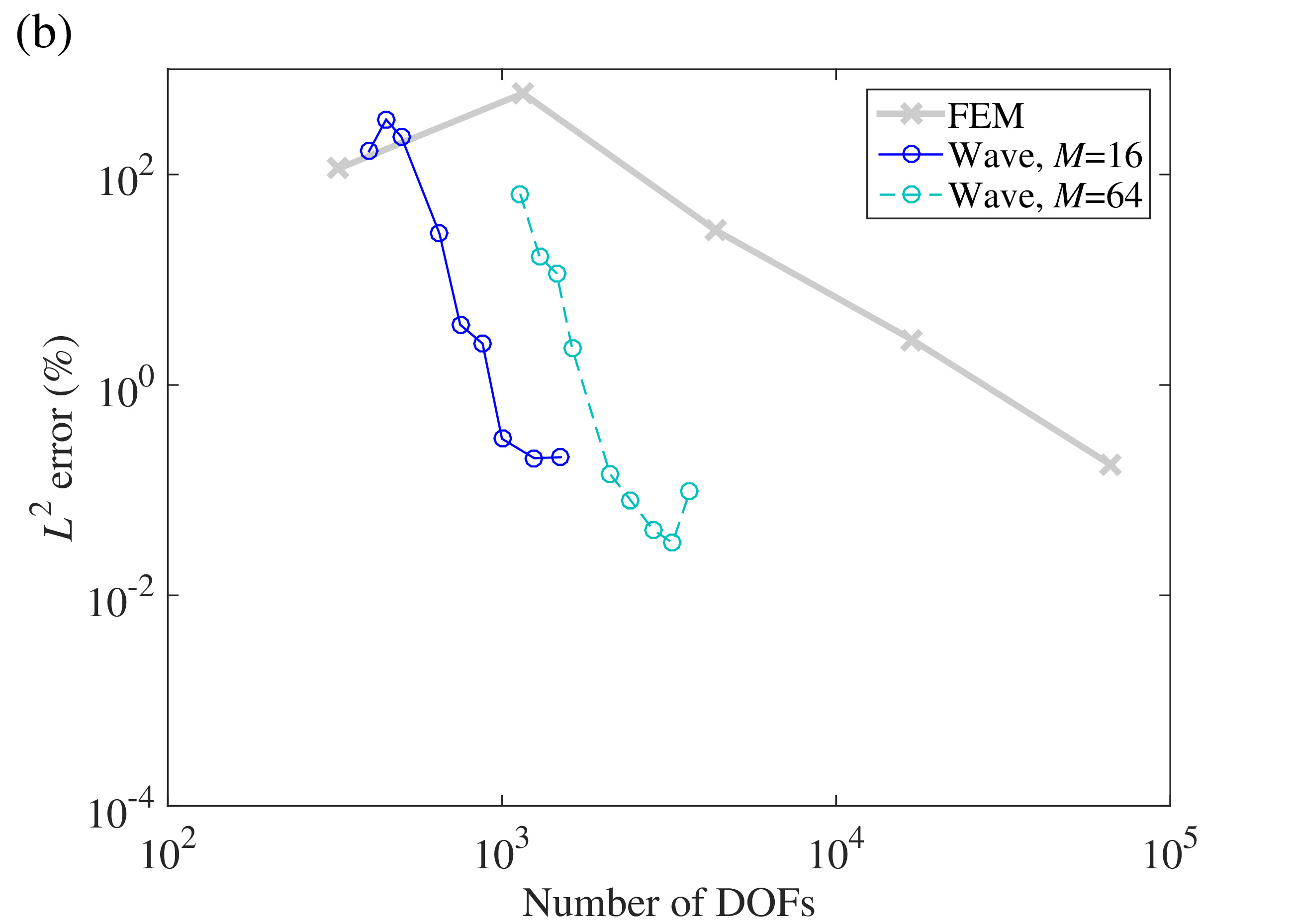}
		\end{center}
		\caption{Convergence curves of classical FEM with mesh refinement and PUFEM with hybrid enrichment combining different numbers of plane waves ($q$) and cubic polynomials ($p=3$) on meshes with $M=16$ and $M=64$ elements, at (a) 1000Hz and (b) 3500Hz.}
		\label{fig_2d_convergence_fc}
	\end{figure}

\begin{figure}[htb!]
		\begin{center}					\includegraphics[width=0.30\textwidth]{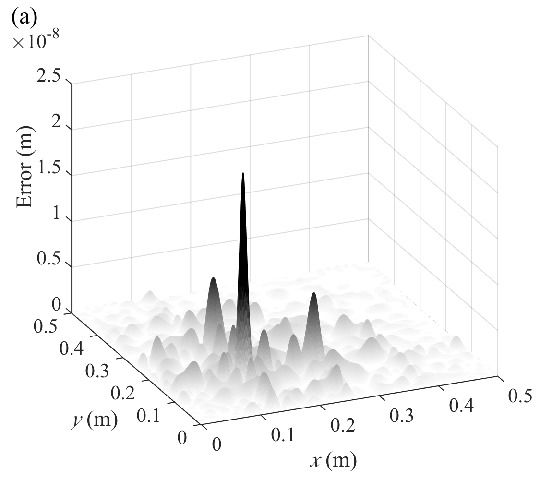}
		\includegraphics[width=0.30\textwidth]{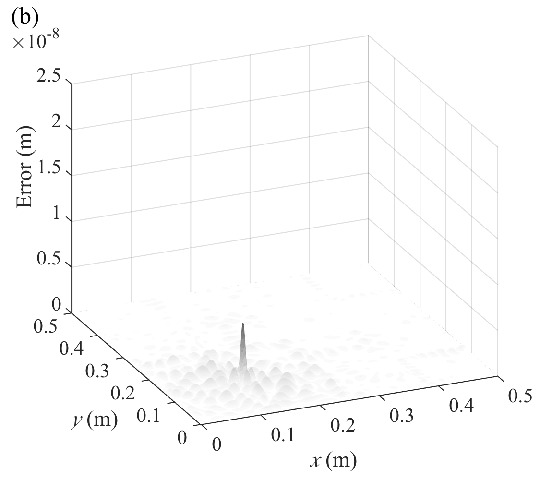}
		\includegraphics[width=0.30\textwidth]{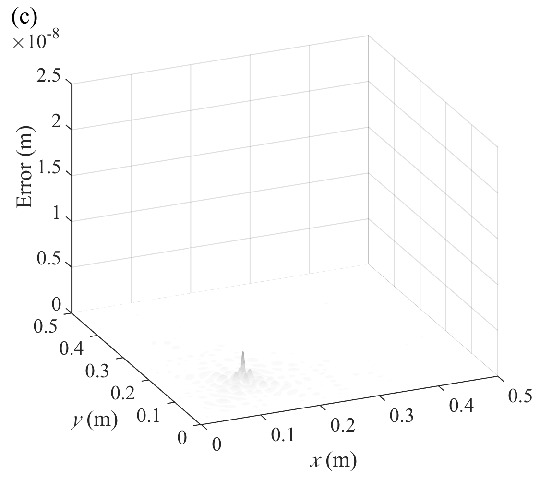}
		\end{center}
		\caption{Absolute displacement difference over the plate at 3500Hz of the PUFEM wave-polynomial enrichment combining $q=40$ with (a) $p=1$ ($\varepsilon_2=0.88\%$), (b) $p=3$ ($\varepsilon_2=0.20\%$) and (c) $p=5$ ($\varepsilon_2=0.089\%$) with $M=16$ elements. }
		\label{fig_2d_error_dist}
	\end{figure}

	Error analyses are conducted on the uniform mesh.
	Figure \ref{fig_2d_convergence_fc} compares the convergence of classical FEM  and PUFEM with hybrid enrichment at two specific frequencies, 1000Hz and 3500Hz. Here, the computational accuracy is estimated via the relative $L^2$ norm error in the two-dimensional domain, defined as:
	\begin{equation}\label{L2_error}
	\varepsilon_2 = \dfrac{\sqrt{ \int_{\Omega} \vert W_{computed} - W_{ref} \vert^2 \mathrm{d}x\mathrm{d}y} }{\sqrt{ \int_{\Omega} \vert W_{ref} \vert^2 \mathrm{d}x\mathrm{d}y } } \times 100 \% \; .
	\end{equation}
	where $W_{ref}$ is the reference solution computed analytically using modal superposition.
	The $L^2$ errors are plotted against the number of DoFs for displacement interpolation.
	Meshes are refined consequently subdividing one primary square element with four sub square elements of equal size and two PUFEM meshes with $M=4\times4$ and $M=8\times8$ elements are investigated in order to evaluate the effect of the PUFEM element size on the performance of the method.
    Results bear similarities with the 1D case except that convergence curves are more steep which is in line with a $q$-refinement, i.e. the number of plane waves increases while the PUFEM mesh is fixed, as opposed to the $h$-refinement shown in Fig. \ref{fig_1d_convergence} where the element size is reduced while the number of DOFs per node is fixed. Two more observations can be made: first, the best results in terms of data reduction is obtained for large size elements, in agreement with previous studies; and second, PUFEM convergence curves all reach a plateau which is attributed to the behavior of the displacement field in the vicinity of the load point. In order to identify this more precisely, it is instructive to remind the analytic solution of an infinite plate subject to a point harmonic \cite{fahy2007}:
	\begin{equation}
	W=F[H_0^{(2)}(kr)-H_0^{(2)}(-jkr)]/8jDk^2, 
	\end{equation}
	where $H_0^{(2)}$ is the Hankel function of the second kind and $r$ is the radial distance from the point force. Thus, the displacement is continuous and behaves locally as $W \approx A + B r^2 \ln(r)$. This means that the second spatial derivatives of the displacement field exhibit a weak singular behavior (of logarithmic type) which is sufficient to affect slightly the convergence curves, however the field remains sufficiently regular to be approximated with a PUFEM wave expansion basis as opposed to the simulation of acoustic pressure fields 
    which exhibit a much stronger singularity in the vicinity of a point source. This has been shown to limit the efficiency of the PUFEM and the problem can be circumvented by considering the scattered field instead (see Ref. \cite{YANG2018} for more details).
    Figure~\ref{fig_2d_error_dist} shows that increasing the polynomial order $p$ while keeping the number plane waves fixed, here $q=40$, can reduce the error locally. Such accuracy cannot be reached by increasing the number of wave directions. 

    \begin{figure}[htb!]
		\begin{center}	\includegraphics[width=0.48\textwidth]{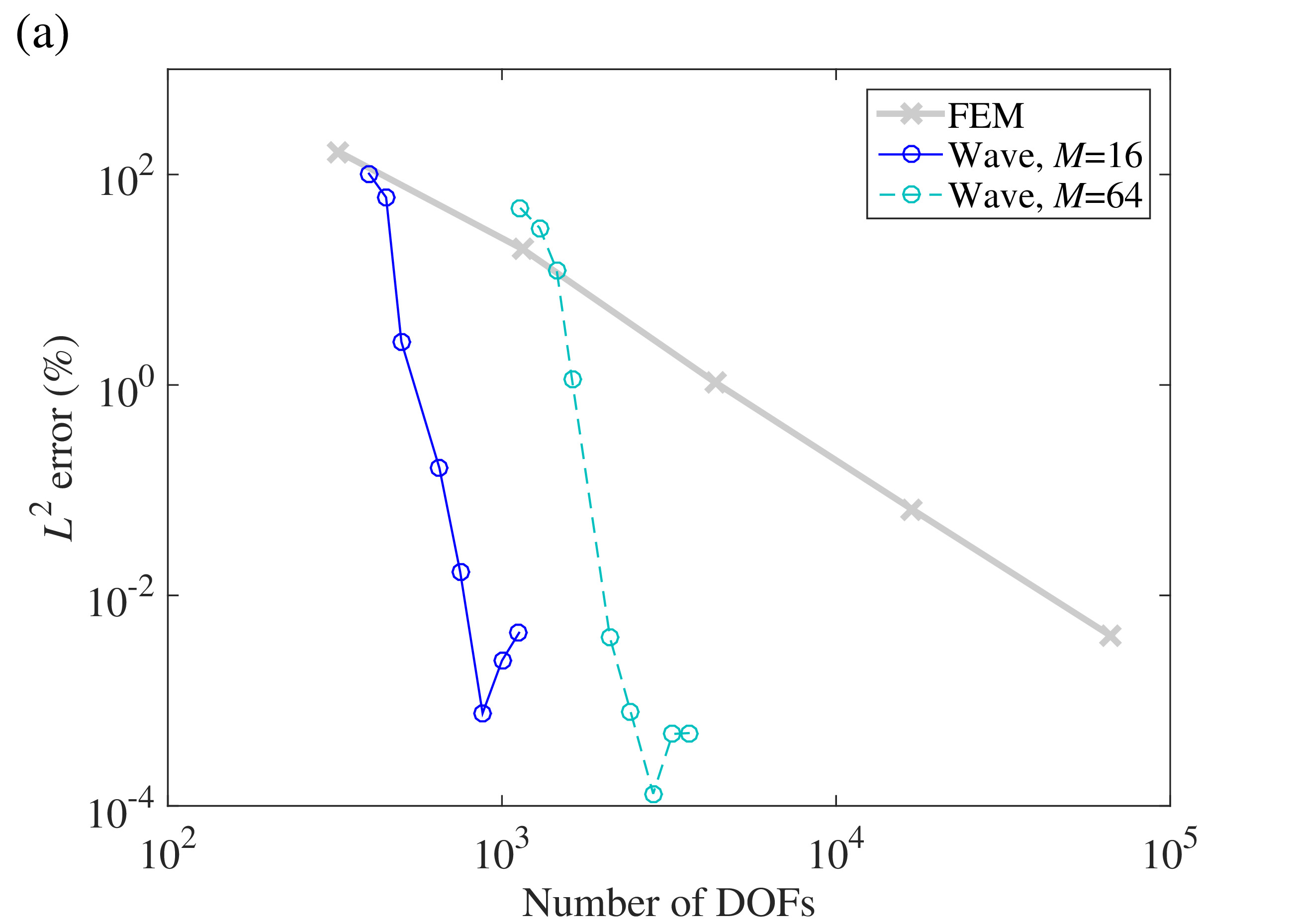}
		\includegraphics[width=0.48\textwidth]{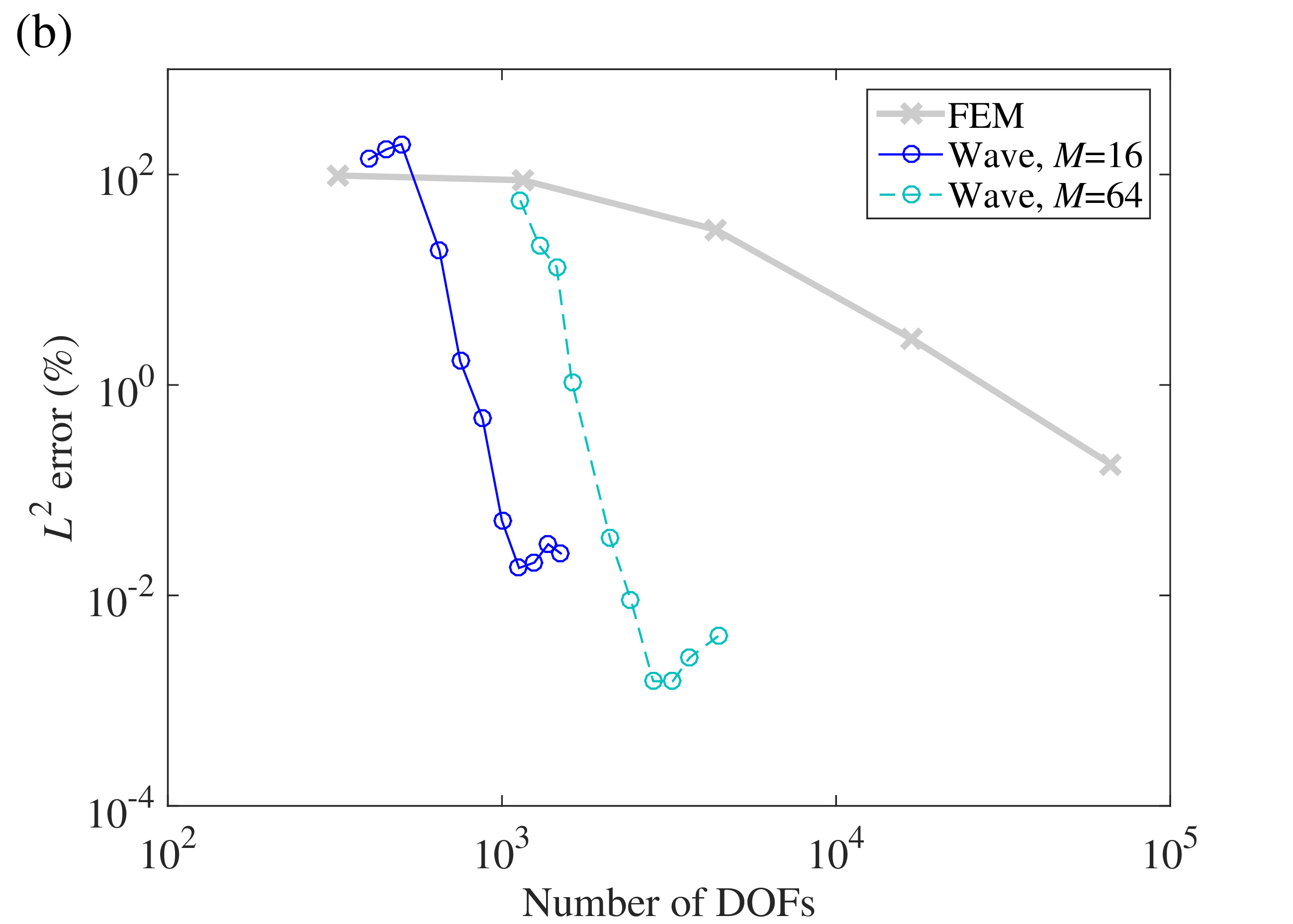}
		\end{center}
		\caption{Convergence curves in the case of a uniform loading at (a) 1000Hz and (b) 3500Hz. Results are obtained by classical FEM with mesh refinement and PUFEM with hybrid enrichment combining different numbers of plane waves ($q$) and cubic polynomials ($p=3$) on two different PUFEM meshes with $M=16$ and $M=64$ elements.}
		\label{fig_2d_convergence_fz}
	\end{figure}
	
	\begin{figure}[htb!]
		\begin{center}	\includegraphics[width=0.48\textwidth]{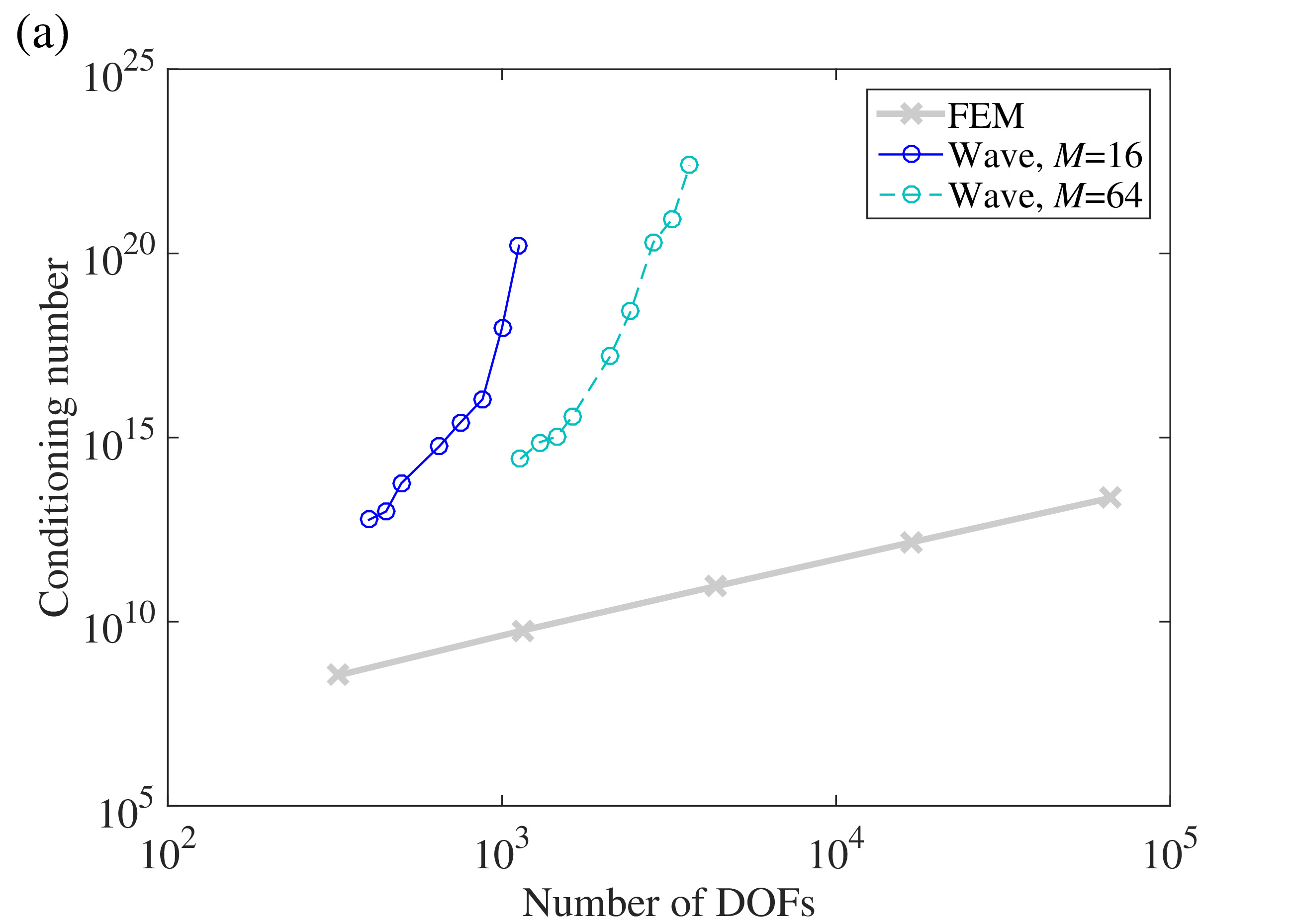}
		\includegraphics[width=0.48\textwidth]{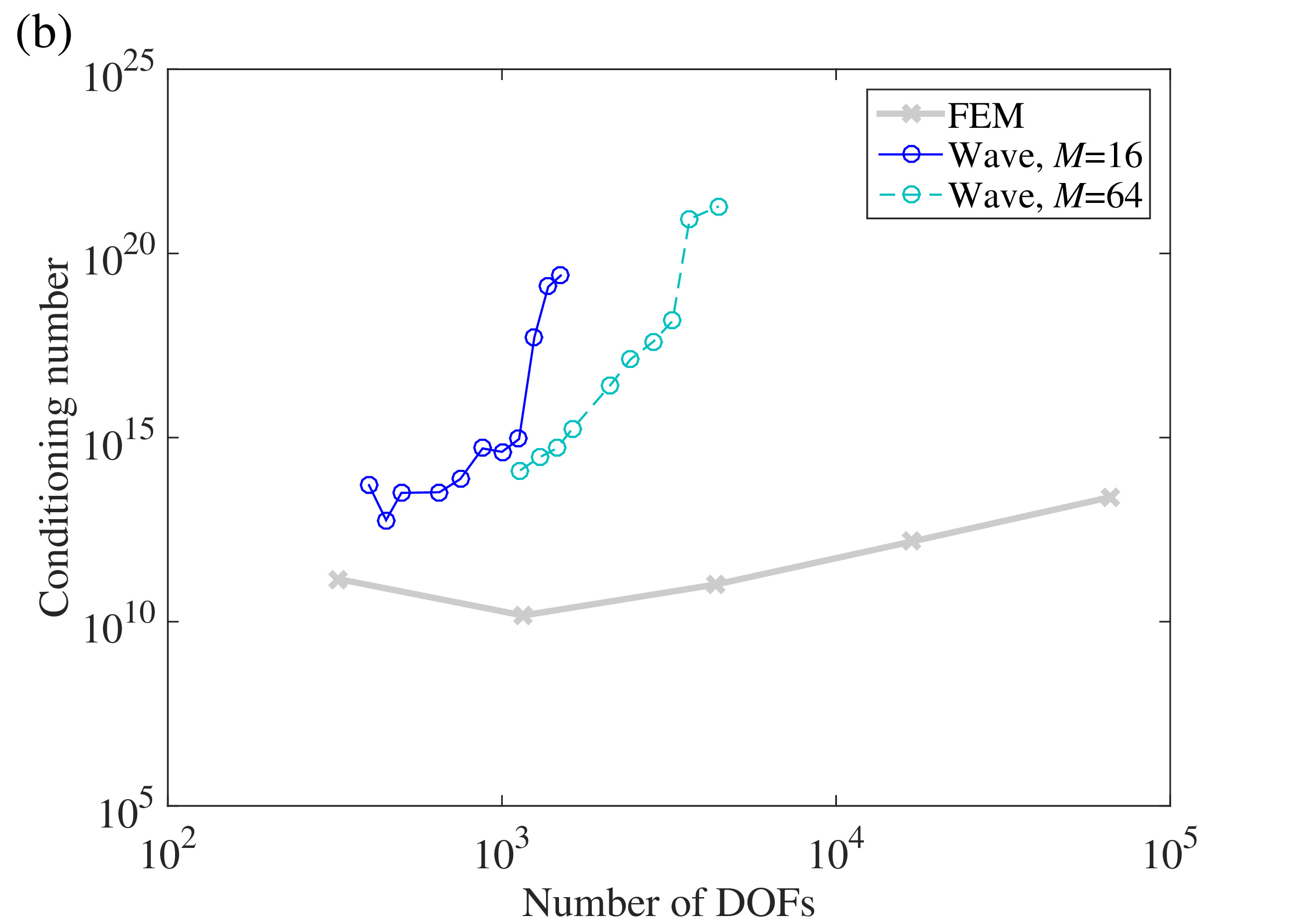}
		\end{center}
		\caption{The conditioning numbers of system matrices $\textbf{K}_{WW}$ at (a) 1000Hz and (b) 3500Hz corresponding to Fig.\ref{fig_2d_convergence_fz}.}
		\label{fig_2d_conditioning_fz}
	\end{figure}

	Convergence results are now given in Fig.~\ref{fig_2d_convergence_fz} for the case of a uniform loading with an unit amplitude, $f_z=1$.
	Results show similar trends as those of Fig.~\ref{fig_2d_convergence_fc} except that errors using PUFEM are reduced at least by an order of magnitude. A closer analysis reveals that the accuracy is limited by the ill-conditioning issues inherent to PUFEM with plane waves \cite{Bettess} and numerical convergence is reached as soon as the conditioning number of the algebraic system exceeds $10^{16}$ and this is confirmed in Fig.~\ref{fig_2d_conditioning_fz}. Similar observations were made in Ref. \cite{YANG2018}.
	However, in the present scenario, the presence of the evanescent waves stemming from the edges of the plate have a detrimental effect on convergence and this seems to worsen as frequency increases. One way to remedy this is to refine the PUFEM mesh which inevitably will reduce the efficiency of the method in terms of data reduction. Another strategy is to exploit the hybrid enrichment and increase the polynomial order.
	For instance, numerical simulations show that the lowest error, here around  $0.02\%$, reached by the wave-polynomial with $p=3$ on $M=16$ PUFEM elements at 3500Hz can be further reduced below $0.01\%$  by simply taking $p=5$. On the contrary, choosing $p=1$ in the formulation does not allow to reach less than $0.6\%$.
	

	
	In order to see this more clearly, two performance indices are introduced to further demonstrate the computational performance of the PUFEM wave-polynomial enrichment for the simulation of bending waves.
	Given a uniform mesh made with square elements, we can define the number of wavelengths spanning over a characteristic element length $h=h_x=h_y$ as 
	\begin{equation}
	\kappa = h/\lambda_b,
	\end{equation}
	By calling $N_{\mathrm{dof}}$ the total number of DOFs needed to interpolate the displacement, the average discretization level $\tau$ that measures the number of DOFs used for modeling a single wavelength \cite{Laghrouche2, Chazot2013_pufem1} is defined as
	\begin{equation}
	\tau=\lambda_b\sqrt{N_{\mathrm{dof}}/S},
	\end{equation}
	where $S=L\times L$ is the surface area of the plate.
	Consider results of Fig.~\ref{fig_2d_point}(a) for instance, at the highest frequency $4000$ Hz, we have  $\kappa=1.8$ and $\tau=4.4$. This value must be compared to the classical rule of thumb requirement that at least `ten nodal points per wavelength' should be used with conventional FEM.
	A $\kappa$-$\tau$-$\varepsilon_2$ convergence analysis is carried out for several frequencies identified by the non-dimensional wavenumbers $kh$ and results are reported in Table \ref{tab_efficiency}.
	The plate is subject to a uniform loading and discretized with 16 elements.
	It can be seen that the PUFEM hybrid enrichment can capture multiple wavelengths per element by using a small number of variables, which is a typical feature of the enrichment by wave solutions. 
	The data reduction, characterized by a low value of $\tau$, is more significant when the wavelength becomes shorter or more precisely when each element contains more wavelengths and at very high frequency $kh=30$, less than 3 DOFs per wavelengths are needed to cap the errors below $1\%$.
	Results of Table \ref{tab_efficiency} also show that errors quickly reach numerical convergence as $q$ increases
	and the associated errors tend to grow for large values of $kh$. 
	A plausible reason for this is the fact that if the finite superposition of propagating plane waves as defined in Eq. (\ref{wave_angle}) can in principle simulate evanescent fields (see a recent paper \cite{Christophe2020} in the case of duct acoustics and also theoretical results in Ref. \cite{perrey2006}), other solutions of the dispersion equation which are of a different nature 
    are deliberately not included in the PUFEM enrichment.
    However, improvements can be made by simply increasing the polynomial order $p$ whilst maintaining 
    very good performance in terms of data reduction.
    From these results, it transpires that (i) the number of wave directions is proportional to the non-dimensional frequencies $kh$, here we have roughly $q \approx 2kh$ to get $\varepsilon_2 \approx 1\%$, which is in line with the acoustic case \cite{Chazot2013_pufem1} and (ii) a two-dimensional polynomial of low-to-moderate degree in order to account for evanescent terms.


\begin{table}	
 \begin{minipage}{0.45\textwidth}
 \small
  \centering
		\begin{tabular}{c c c c c c}
		    \hline
			$kh$ & $p$ & $q$ & $\kappa$ & $\tau$ & $\varepsilon_2(\%)$ \\\hline
			 5 & 3 & 15 & 0.80 & 7.85 & 0.17  \\
			    &   & 20 & 0.80 & 8.60 & 0.011 \\
			    &   & 25 & 0.80 & 9.29 & 0.0032 \\\hline			
			 10 & 3 & 25 & 1.59 & 4.65 & 0.56  \\
			    &   & 30 & 1.59 & 4.97 & 0.036 \\
			    &   & 35 & 1.59 & 5.27 & 0.097 \\\hline
			 15 & 3 & 30 & 2.39 & 3.31 & 4.47 \\
			    &   & 40 & 2.39 & 3.70 & 0.074 \\
			    &   & 50 & 2.39 & 4.06 & 0.049 \\\hline
			 20 & 3 & 30 & 3.18 & 2.48 & 125.05 \\
			    &   & 45 & 3.18 & 2.91 & 2.20 \\
			    &   & 60 & 3.18 & 3.29 & 2.12 \\\hline	
		\end{tabular}
  \end{minipage}		
  \begin{minipage}{0.45\textwidth}
   \centering
   \small   
		\begin{tabular}{c c c c c c}
		\hline
		$kh$ & $p$ & $q$ & $\kappa$ & $\tau$ & $\varepsilon_2(\%)$ \\\hline
			 5 & 5 & 15 & 0.80 & 9.42 & 0.015  \\
			    &   & 20 & 0.80 & 10.06 & 0.0014 \\
			    &   & 25 & 0.80 & 10.65 & 0.0024 \\\hline			
			 10 & 5 & 25 & 1.59 & 5.33 & 0.085  \\
			    &   & 30 & 1.59 & 5.61 & 0.036 \\
			    &   & 35 & 1.59 & 5.88 & 0.0053 \\\hline		
			 15 & 5 & 30 & 2.39 & 3.74 & 0.98 \\
			    &   & 40 & 2.39 & 4.09 & 0.0065 \\
			    &   & 50 & 2.39 & 4.41 & 0.0039 \\\hline
			 20 & 5 & 30 & 3.18 & 2.80 & 39.91 \\
			    &   & 45 & 3.18 & 3.19 & 0.19 \\
			    &   & 60 & 3.18 & 3.53 & 0.28 \\\hline
			    
		\end{tabular}
	   \end{minipage}
	   \vspace{0.2cm}
	   
 \begin{minipage}{0.45\textwidth}
 \small
  \centering
		\begin{tabular}{c c c c c c}
		    \hline
			$kh$ & $p$ & $q$ & $\kappa$ & $\tau$ & $\varepsilon_2(\%)$ \\\hline
			 25 & 7 & 40 & 3.98 & 2.74 & 11.77  \\
			    &   & 50 & 3.98 & 2.91 & 0.40 \\
			    &   & 60 & 3.98 & 3.08 & 0.40 \\
			    &   & 70 & 3.98 & 3.23 & 0.20	\\\hline			
			 30 & 7 & 50 & 4.77 & 2.43 & 17.64  \\
			    &   & 60 & 4.77 & 2.57 & 2.85 \\
			    &   & 70 & 4.77 & 2.70 & 3.45 \\ 
			    &   & 80 & 4.77 & 2.82 & 1.72 \\\hline
		\end{tabular}
  \end{minipage}
   \begin{minipage}{0.45\textwidth}
 \small
  \centering
		\begin{tabular}{c c c c c c}
		    \hline
			$kh$ & $p$ & $q$ & $\kappa$ & $\tau$ & $\varepsilon_2(\%)$ \\\hline
			 25 & 9 & 40 & 3.98 & 3.06 & 5.24  \\
			    &   & 50 & 3.98 & 3.22 & 0.28 \\
			    &   & 60 & 3.98 & 3.37 & 0.077 \\
			    &   & 70 & 3.98 & 3.51 & 0.068	\\\hline			
			 30 & 9 & 50 & 4.77 & 2.68 & 2.02  \\
			    &   & 60 & 4.77 & 2.81 & 0.57 \\
			    &   & 70 & 4.77 & 2.93 & 0.49 \\ 
			    &   & 80 & 4.77 & 3.04 & 0.43  \\\hline
		\end{tabular}
  \end{minipage}


\caption{Performances of the PUFEM with wave-polynomial enrichment for the simply supported square plate meshed with $M=16$ elements.}
\label{tab_efficiency}
\end{table}

Further analyses are carried out to illustrate the role of the added polynomials in the hybrid enrichment especially in the vicinity of the edges of the plate where evanescent waves are expected to dominate.
The order of polynomials associated with nodes inside the plate ($p_i$) is chosen to be lower than those associated with nodes lying on the edges ($p_e$). 
Numerical tests were carried out for two frequencies $kh=20, 30$ with results tabulated in Table \ref{tab_adaptive}. 
It is observed that similar accuracy can be obtained while  reducing the polynomial enrichment for internal nodes.
This shows the interest of using an adaptive polynomial enrichment scheme which  leads to a smaller number of degrees of freedom and a better averaged discretization level.

\begin{table}

 \small
  \centering
		\begin{tabular}{c c c c c c c}
		    \hline
			$kh$ & $p_e$ & $p_i$ & $q$ & $\kappa$ & $\tau$ & $\varepsilon_2(\%)$ \\\hline
			 20 & 5 & 3 & 30 & 3.18 & 2.69 & 43.92 \\
			    &   &   & 45 & 3.18 & 3.09 & 0.98 \\
			    &   &   & 60 & 3.18 & 3.45 & 0.35 \\\hline

			 20 & 5 & 1 & 30 & 3.18 & 2.62 & 63.85 \\
			    &   &   & 45 & 3.18 & 3.03 & 2.44 \\
			    &   &   & 60 & 3.18 & 3.39 & 0.25 \\\hline

			 30 & 9 & 7 & 50 & 4.77 & 2.59 & 7.34 \\
			    &   &   & 60 & 4.77 & 2.72 & 0.77 \\
			    &   &   & 70 & 4.77 & 2.85 & 0.68 \\ 
			    &   &   & 80 & 4.77 & 2.96 & 0.37 \\\hline			    
			 30 & 9 & 5 & 50 & 4.77 & 2.52 & 4.65  \\
			    &   &   & 60 & 4.77 & 2.65 & 1.34 \\
			    &   &   & 70 & 4.77 & 2.78 & 0.46 \\ 
			    &   &   & 80 & 4.77 & 2.90 & 0.48  \\\hline			    
		\end{tabular}
	
\caption{Performance of the hybrid enrichment using reduced-order polynomials at internal nodes of the simply supported square plate meshed with $M=16$ elements ($p_e$ and $p_i$ stand for enrichment orders for edge and internal nodes, respectively).}
\label{tab_adaptive}
\end{table}

\subsection{Application to an L-shaped plate}	

\begin{figure}[htb!]
		\begin{center}	\includegraphics[width=0.48\textwidth]{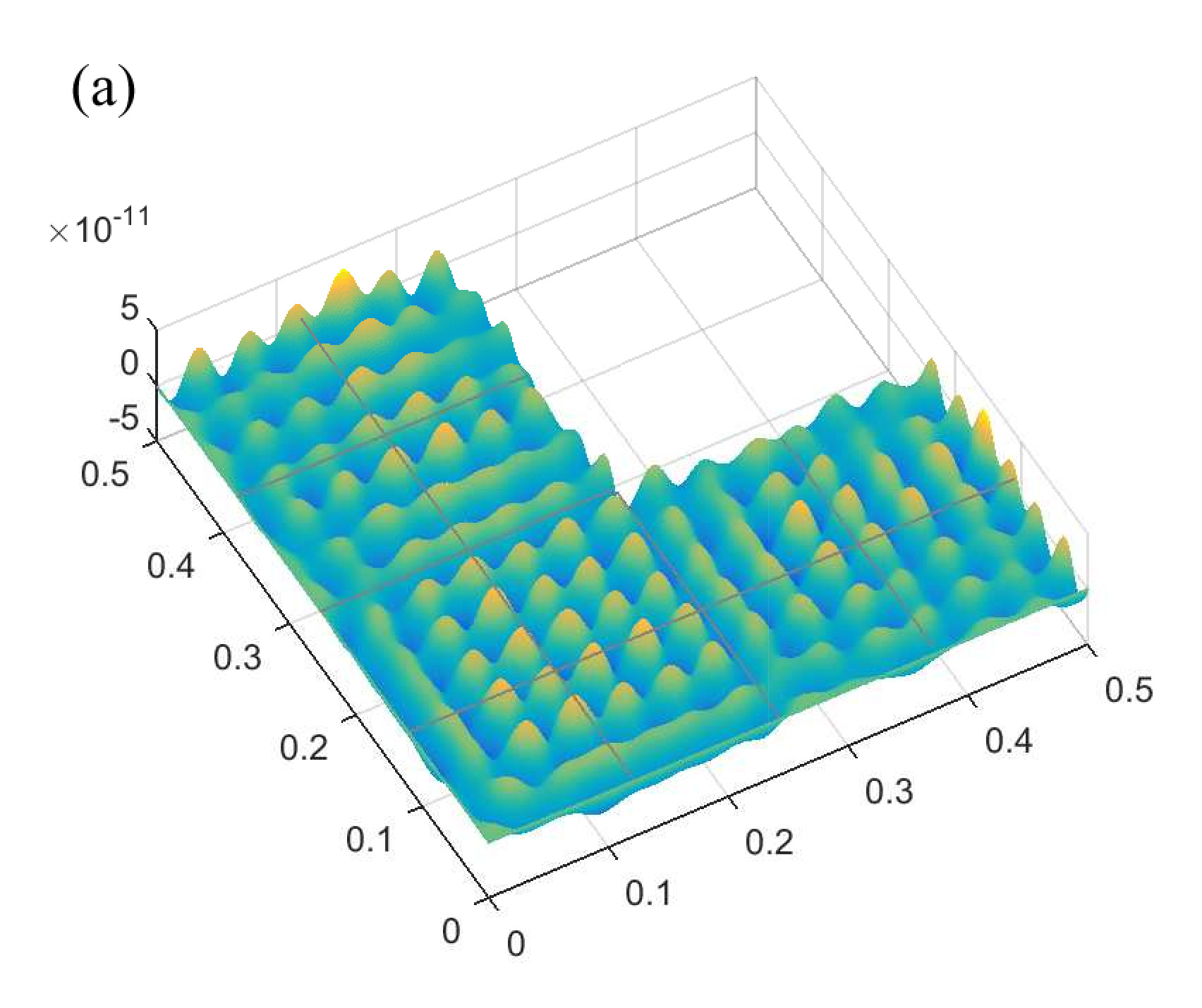}
		\includegraphics[width=0.48\textwidth]{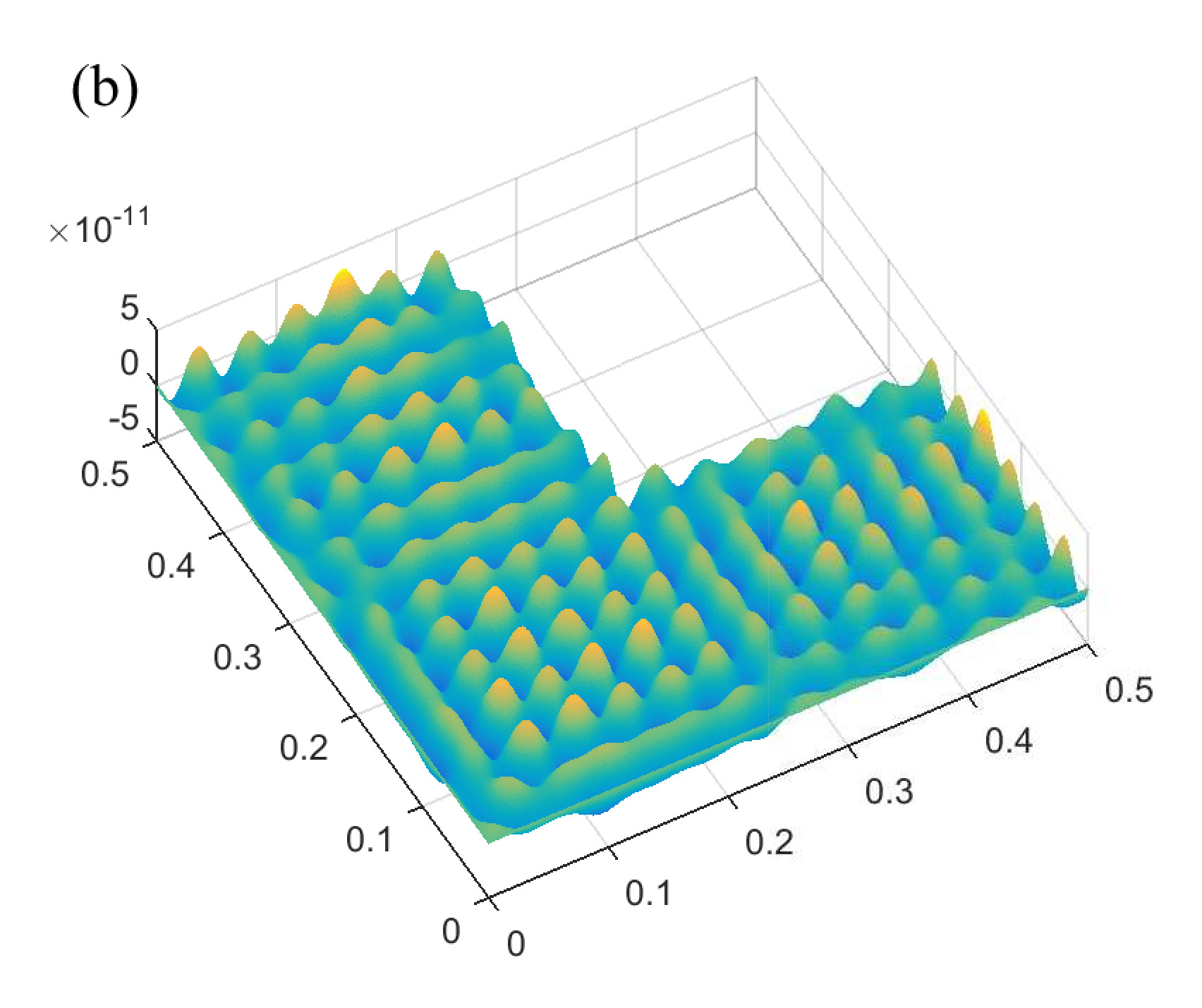}
		\end{center}
		\caption{Responses of an L-shaped plate calculated by PUFEM with wave-polynomial enrichment and FEM with a highly refined mesh ($kh=20$). }
		\label{fig_L_shape_plate}
	\end{figure}

The present numerical scheme is applied to a non-square plate. Here, we consider an L-shaped plate by removing a quarter of the previous square plate ($x>L/2$ and $y>L/2$). The plate is subject to a uniformly distributed loading with $f_z=1$. Two longest edges are simply supported while others are free of constrains.
The plate is meshed by square elements of the same size.
As shown in Fig.~\ref{fig_L_shape_plate}, there is a good agreement of the responses calculated by the PUFEM with wave-polynomial enrichment ($q=60, p=7$) using 12 PUFEM elements and the classical FEM using a highly refined mesh with 49152 CR elements. Here, $kh=20$ and
the numbers of DoFs are $2016$ and $198660$, respectively.
The relative difference of the responses is around $0.5\%$, evaluated by Eq. (\ref{L2_error}) with the FEM solution as the reference.
Clearly, PUFEM outperforms classical FEM with CR elements since comparable accuracy can be obtained with nearly one hundred times less DoFs and at higher frequency, the gain in terms of data reduction remains similar.

	\section{Conclusions}
	
	In this paper, we have presented a conforming PUFEM thin plate bending element for the modeling of steady-state structural harmonic vibrations.
	Hermite shape functions associated with the displacement are used to form a partition of unity which, by construction, guarantees that first derivatives  are continuous everywhere.
	The method allows to consider a variety of enrichment functions. In the present work, 
	high-order polynomials and plane wave solutions to the differential equations are considered. 
	This leads to a hybrid wave-polynomial enrichment which is developed and applied 
	in both 1D and 2D configurations.
	The one-dimensional analysis permits to illustrate and assess the performance of the method. It is found that polynomial enrichment with sufficiently high order yields better convergence than classical FEM elements, and the hybrid enrichment which includes two waves propagating in opposite directions leads to the best computational performance in terms of both accuracy and data reduction.



    The PUFEM with hybrid enrichment is extended to the simulation of flexural vibrations in plates by combining plane waves propagating in various directions with additional polynomial terms.
    Specific constraints which are prescribed for the edges of the plate are imposed using Lagrange multipliers. The method is applied and tested using simply supported rectangular plates with singular and uniformly distributed loads.
    It is found that the method can achieve accurate predictions over a wide frequency range by using a relatively small number of basis functions, when compared to classical CR elements, and has the capacity to accommodate certain mesh irregularity. 
    A convergence analysis carried out in comparison with classical CR elements shows that the PUFEM exhibits extremely high rate of convergence. It is also found that the polynomial enrichment allows capturing evanescent or nearly-singular fields in the vicinity of a singular load point and along the edges of the plates. 
    
    It is shown that the Partion of Unity formulation presented in this work permits to construct a wide family of element shape functions with sufficient regularity for solving bending waves in beams and plates at high frequencies and there are good reasons to believe that other type of functions could also be tailored for more specific configurations involving other geometrical features for instance.
    Our future works include the application and the extension of thin plate PUFEM elements to take into account variable structural thickness \cite{ma2018, huang2016, zhao2019, pelat2020}.
	


	\section{References}
	
	\bibliographystyle{unsrt}
	\bibliography{beam_pufem}

\end{document}